\documentclass{article}

\usepackage{arxiv}

\usepackage[utf8]{inputenc} 
\usepackage[T1]{fontenc}    
\usepackage{hyperref}       
\usepackage{url}            
\usepackage{booktabs}       
\usepackage{amsfonts}       
\usepackage{nicefrac}       
\usepackage{microtype}      

\usepackage{graphicx}

\usepackage{colortbl}		
\usepackage{multicol}        

\usepackage{amsbsy}
\usepackage{amsfonts}
\usepackage{amsmath}

\usepackage{natbib,har2nat}

\usepackage{algorithm}
\usepackage[algo2e, boxed, ruled, vlined]{algorithm2e}

\usepackage{mathabx}
\usepackage{booktabs}       
\usepackage{multirow}

\usepackage{color}
\usepackage{graphics,hyperref}
\usepackage{latexsym}
\usepackage{times, bm}

\usepackage{amssymb,amsmath,mathrsfs,amsfonts,amsthm, amsbsy}

\newtheorem{corollary}{Corollary}[section]   
\newtheorem{remark}{Remark}  
\newtheorem{lemma}{Lemma}[section]

\newtheorem{proposition}{Proposition}[section]
\newtheorem{condition}{Condition}

\usepackage{color}

\DeclareMathOperator*{\argmax}{arg\,max}

\newcommand{\bigO}{\ensuremath{\mathcal{O}}}
\newcommand{\bhat}[1]{\expandafter\hat#1}

\newcommand{\bo}[1]{\mathbf{#1}}

\newcommand{\R}{\mathbb{R}} 

\newcommand{\z}{\bo z} 
\newcommand{\x}{\bo x}
\newcommand{\X}{\bo X}
\newcommand{\y}{\bo y}
\newcommand{\G}{\bo \Gamma} 
\newcommand{\g}{\boldsymbol{\gamma}} 
\newcommand{\Om}{\bo \Omega} 
\newcommand{\I}{\bo I} 
\newcommand{\U}{\bo U}
\newcommand{\bu}{\bo u}
\newcommand{\bv}{\bo v}
\newcommand{\V}{\bo V}
\newcommand{\Sig}{\bo \Sigma}

\newcommand{\T}{\bo T}
\newcommand{\D}{\bo D}
\newcommand{\s}{\bo s}
\newcommand{\h}{\bo h}
\newcommand{\e}{\bo e}
\newcommand{\0}{\bo 0}
\newcommand{\M}{\bo M}

\def\Cov{\mathop{\mathrm{Cov}}\nolimits}

\def\n{\frac{1}{n}\sum_{i=1}^n}

\def\wM{\widehat{\bo M}}
\def\wG{\widehat{\bo \Gamma}}

\def\wD{\widehat{\bo D}}

\title{Test of the Latent Dimension of a Spatial Blind Source Separation Model}

\date{} 					

\author{Christoph~Muehlmann \\
	Institute of Statistics \& Mathematical Methods in Economics \\
	Vienna University of Technology, Austria \\
	\texttt{christoph.muehlmann@tuwien.ac.at} \\
	\And
	Fran\c{c}ois~Bachoc\\
    Institut de Math\'{e}matiques de Toulouse \\
    Universit\'{e} Paul Sabatier, France \\
    \texttt{francois.bachoc@math.univ-toulouse.fr}
	\And
	Klaus~Nordhausen \\
	Institute of Statistics \& Mathematical Methods in Economics \\
	Vienna University of Technology, Austria \\
	\texttt{klaus.nordhausen@tuwien.ac.at} \\
	\And
    Mengxi~Yi \\
    School of Statistics \\
    University of International Business and Economics, China \\
	\texttt{myi@uibe.edu.cn} \\}

\renewcommand{\headeright}{}
\renewcommand{\undertitle}{}

\begin{document}
\maketitle

\begin{abstract}
We assume a spatial blind source separation model in which the observed multivariate spatial data is a linear mixture of latent spatially uncorrelated Gaussian random fields containing a number of pure white noise components. 
We propose a test on the number of white noise components and obtain the asymptotic distribution of its statistic for a general domain. We also demonstrate how computations can be facilitated in the case of gridded observation locations. Based on this test, we obtain a consistent estimator of the true dimension. Simulation studies and an environmental application demonstrate that our test is at least comparable to and often outperforms  bootstrap-based techniques, which are also introduced in this paper.
\end{abstract}

\section{Introduction} \label{Sec:Intro}
With the advance of technology, massive amounts of multivariate spatial data can be collected. As one example, researchers may use these datasets to investigate various issues in geographical, ecological (\cite{legendre2012numerical}) or atmospheric (\cite{von2001statistical}) sciences. Commonly, in these spatially correlated datasets there are dependencies both within and among the individual data processes, which makes modeling the multivariate spatial data a challenge. The difficulty of modeling is further intensified when the dimensionality $p$ is large. With a dataset of size $n$, it takes a total of $cp(p+1)/2$ parameters to describe the full covariance and cross-covariance structure of the model, where $c$ is the number of characteristic parameters per covariance and cross covariance. Further, it requires a computational cost of $\bigO(n^3p^3)$ for prediction using optimal linear predictors and for Gaussian likelihood evaluation
(\cite{cressie2015statistics,legendre2012numerical}).

One way to approach this problem is to use the spatial blind source separation (SBSS) framework see \cite{nordhausen2015blind} and \cite{bachoc2018spatial}. Blind source separation (BSS) is a well-studied multivariate procedure used to recover latent variables when only a linear mixture of them is observed; see for example \cite{comon2010handbook,NordhausenOja2018}. 
 A common assumption for BSS is that the latent variables are second-order stationary and uncorrelated. That is, we assume $\x(\s)=\Om \z (\s)$, where $\x (\s) \in\R^p$ is the observed $p$-variate measurement at location $\s \in \mathbb{R}^d$, $\z (\s) \in\R^p$ is a latent second-order stationary $p$-variate source with uncorrelated components, and $\Om\in\R^{p\times p}$ is an unknown full-rank mixing matrix. To estimate the unmixing matrix $\G$, i.e. $\Om^{-1}$, \cite{nordhausen2015blind} have proposed an estimator based on the simultaneous diagonalization of two scatter matrices, and \cite{bachoc2018spatial} have extended this method to jointly diagonalize more than two scatter matrices for multivariate spatial data. 

The SBSS model of \cite{nordhausen2015blind} gives no preference to any of the latent components, with all $p$ of them being basically of equal interest. However, in practical cases of BSS, it is often assumed that only a few components are of interest and to be regarded as the signal, while the remaining components are discarded as noise. That is, it is supposed that the latent components consist of two parts, $\z=(\z_{s}^T, \z_{w}^T)^T$, where $\z_{s}\in\R^q$ is the signal, and $\z_{w}\in\R^{p-q}$ is the noise. \cite{bootstrap_ts, virta2018determining} and \cite{NordhausenVirta2018} all consider components with serial dependence as signals in a time series context.

Following the work of \cite{virta2018determining}, in this study we consider SBSS in which signals are characterized  as components  having second-order spatial dependence.
We derive a test for the signal dimension $q$ based on the joint diagonalization of two or more scatter matrices that are specified by kernel functions, as done by \cite{bachoc2018spatial}. We then provide the asymptotic distribution of the test statistic. This asymptotic result does not rely on the strong assumption made by \cite{bachoc2018spatial}, that all the signal and noise components are asymptotically identifiable. To avoid this assumption, new proof techniques are needed, and we rely in particular on the extension of arguments made by \cite{virta2018determining} to a spatial setting. In addition, we demonstrate that using a different normalization factor for the scatter matrices than that used by \cite{bachoc2018spatial} enables the obtainment of a neater asymptotic distribution of the test statistic (see Remark \ref{remark:renormalization}). Based on the test, we then provide a consistent estimator of the unknown signal dimension.

We put forward several bootstrap versions of the test. For both the asymptotic and bootstrap tests, we demonstrate that computational gains are obtained when the observation locations are gridded. In an extensive simulation study, we then show that the various tests already have levels close to the nominal one, for small to moderate sample sizes. We also observe an accurate estimation of the signal dimension. We conclude that the asymptotic test is comparable to and often outperforms the bootstrap ones while being less computationally demanding. Employing an environmental application, we then show that our methods enable the reduction of the dimension of a multivariate spatial data set, retaining the most interpretable and informative estimated independent components and discarding the unusable ones as noise.

The remainder of the paper is organized as follows. In Section \ref{Sec:Setup}, we introduce the statistical setting of the problem and present our test statistic. The methods and main results are then described in Section \ref{Sec:Theory}, the simulation results are reported in Section \ref{Sec:Sim}, and the environmental application is outlined in Section \ref{Sec:Ex}. We finally discuss some concluding remarks in Section \ref{Sec:Conclusion}. The proofs of the theoretical results are provided in the appendix.

\section{Setup and Model}\label{Sec:Setup}

Suppose our data consists of a $p$-dimensional multivariate random field $\x(\s)=\{x_1(\s), \cdots, x_p(\s)\}^T$, $\s \in \mathcal{S}$, where $\mathcal{S}\subseteq\R^d$ is a region of interest. The covariance and cross-covariance functions of 
$\x$, defining its second-order structure, are some of its central characteristics. We can refer, for instance, to \cite{deIaco2013using,genton2015cross} and \cite{gneiting2010matern} for an introduction and various approaches to modeling these.

Here, the second-order structure of $\x$ is assumed to obey an SBSS model:
\begin{equation}\label{eq:model}
\x(\s)=\Om \z(\s),
\end{equation}
where $\Om$ is a $p\times p$ unknown invertible matrix, and  $\z(\s)=\{z_1(\s),\cdots,z_p(\s)\}^T$ is the latent field having independent components with $\Cov(\z(\s))=\I_p$ for all $\s\in\mathcal{S}$.

Let $I(\cdot)$ denote the indicator function throughout this paper and consider the kernel functions $f_0, f_1, \cdots, f_k$, with $f_\ell: \R^d\to\R$ for $\ell=0, \cdots, k$, and with $f_0(\s)=I(\s=\boldsymbol{0})$. For $f\in\{f_0, f_1, \cdots, f_k\}$, let
\[
F_{n,f}=\frac{1}{n}\sum_{i,j=1}^nf^2(\s_i-\s_j),
\]
where $\{ \s_1,\ldots,\s_n \} \subseteq \mathcal{S}$ is the set of two-by-two distinct observation points.
Notice that $F_{n,f_0}=1$. For $f\in\{f_1, \cdots, f_k\}$. The population local covariance (or scatter) matrices are then defined as, 
\begin{equation}\label{eq:plc}
\M(f)=\frac{1}{n\sqrt{F_{n,f}}}\sum_{i=1}^n\sum_{j=1}^nf(\s_i-\s_j)\mathbb{E}\left(\x(\s_i)\x(\s_j)^T\right) \quad \text{and} \quad \M(f_0)=\n\mathbb{E}\left(\x(\s_i)\x(\s_i)^T\right),
\end{equation}
and the corresponding sample local covariance matrices are defined as
\begin{equation}\label{eq:slc}
\wM(f)=\frac{1}{n\sqrt{F_{n,f}}}\sum_{i=1}^n\sum_{j=1}^nf(\s_i-\s_j)\x(\s_i)\x(\s_j)^T \quad \text{and} \quad \wM(f_0)=\n\x(\s_i)\x(\s_i)^T.
\end{equation}
\begin{remark} \label{remark:renormalization}
The normalizing quantity $nF_{n,f}^{1/2}$ in \eqref{eq:plc} and \eqref{eq:slc} is slightly different from that in \cite{bachoc2018spatial}, where simply $n$ is used. Here, the introduction of $F_{n,f}^{1/2}$ enables us to obtain a simple and elegant asymptotic distribution of the test statistic for the number of noise components (see Proposition \ref{prop:asymt}).
\end{remark}

The $k+1$ sample local covariance matrices
$
\wM(f_0), \wM(f_1), \cdots, \wM(f_k)
$
are used to estimate the unmixing matrix  $\G$ as 
\begin{equation}\label{eq:hg}
\wG\in\argmax_{\substack{\G: \G\wM(f_0)\G^T=\I_p\\
\G\text{ has rows }\g_1^T,\cdots,\g_p^T\\
(\sum_{\ell=1}^k\{\g^T_j\wM(f_{\ell})\g_j\}^2)_{j=1,\cdots,p}\text{ are in descending order}}} \sum_{\ell=1}^k\sum_{j=1}^p\{\g^T_j\wM(f_{\ell})\g_j\}^2.
\end{equation}
The unmixing matrix should ``diagonalize" all $k$ local covariance matrices and we let for $\ell=1,\cdots,k$,
\[
\wD_{\ell}=\wG\wM(f_{\ell})\wG^T\]
where all $\wD_{\ell}$ should be close to a diagonal matrix. Note that for finite data exact diagonalization is usually possible only for $k=1$. Further, by definition,
$
\sum_{\ell=1}^k\wD^2_{\ell,1,1}\ge\cdots\ge\sum_{\ell=1}^k\wD^2_{\ell,p,p}.$
We are now interested in the case in which there are $q$ ``real" continuous random fields in $\z$, while the remaining $p-q$ components are white noise.

For $q\in\{0,\cdots,p-1\}$, we are interested in testing the following hypothesis
\begin{center}
$H_{0q}:$ There are exactly $p-q$ white noise processes in $\z$.
\end{center}

This hypothesis is formalized in the following two conditions:
\begin{condition}\label{c:1}
For $a=1, \cdots, p-q$, the covariance function of $z_{q+a}$ is given by 
$$\Cov(z_{q+a}(\bu), z_{q+a}(\bv))=I(\bu-\bv=\boldsymbol{0}).$$
\end{condition}
\begin{condition}\label{c:2}
For $\ell=1,\cdots,k, f_{\ell}$ is symmetric and satisfies $f_{\ell}(\bm{0})=0$. 
For $a=1, \cdots, q$, we have
\[
\liminf_{n\to\infty}\sum_{\ell=1}^k\M(f_{\ell})^2_{a,a}>0.\]
\end{condition} 
Note that in Conditions \ref{c:1} and \ref{c:2}, we assume that the sources are ordered such that the $q$ signal components come first and are followed by the $p-q$ noise components. As the order of the sources is not identifiable, this assumption comes without loss of generality. The fulfillment of Condition \ref{c:2} means that the correlation in the signal fields $z_1,\cdots,z_q$ is sufficient for these signals to be asymptotically separated from the noise fields $z_{q+1}, \cdots, z_{p}$. It should also be noted that we do not need to consider the stronger assumption that the $q$ vectors of the $a$-th diagonal elements in $\M(f_1), \cdots, \M(f_k)$, for $a=1,\cdots,q$, for the signal random fields, are asymptotically distinct (see Assumption 9 in \cite{bachoc2018spatial}) and non-zero. Conditions \ref{c:1} and \ref{c:2} motivate the following block decompositions for $\ell=1,\cdots,k$:
\[
\wM(f_{\ell})=\begin{pmatrix}
\wM(f_{\ell})_{qq} & \wM(f_{\ell})_{q0}\\
\wM(f_{\ell})_{0q} & \wM(f_{\ell})_{00}
\end{pmatrix} \quad \text{ and } \quad 
\wD_{\ell}=\begin{pmatrix}
\wD_{\ell,qq} & \wD_{\ell,q0}\\
\wD_{\ell,0q} & \wD_{\ell,00}
\end{pmatrix}
\]
where the blocks $\wM(f_{\ell})_{qq}$ and $\wD_{\ell,qq}$ have size $q\times q$ and the blocks $\wM(f_{\ell})_{00}$ and $\wD_{\ell,00}$ have dimension $(p-q)\times (p-q)$.\\

Then our test statistic is
\begin{equation}\label{eq:t}
t_q=\frac{n}{2}\sum_{\ell=1}^k||\wD_{\ell,00}||^2,
\end{equation}
where $||\cdot||$ is the Frobenius norm. The test statistic is then expected to be bounded under the null hypothesis, and to diverge when one of $z_{q+1},\ldots,z_{p}$ is spatially correlated. The test will reject the null hypothesis $H_{0q}$ if $t_q$ is larger than a certain threshold, in which case the dataset provides indications that more than $q$ signal components are present.  For a nominal level $\alpha \in (0,1)$, the threshold will be set to the quantile $1-\alpha$ of the asymptotic distribution of Proposition \ref{prop:asymt} or \ref{prop:dw} or Corollary \ref{coro:dw}, depending on the context.

\section{Theory and Methodology} \label{Sec:Theory}
\subsection{Asymptotic Tests for Dimension}
Assume now that $\x$ satisfies Model \eqref{eq:model}. Then, let $q$ denote the true value of the signal dimension (i.e., $H_{0q}$ is true) and consider the limiting distribution of $t_q$. To establish the asymptotic results, we need to introduce a few technical conditions: 

\begin{condition}\label{c:3}
The random fields $z_1,\cdots,z_p$ are independent centered stationary Gaussian random fields.
\end{condition}

The zero mean assumption is replaced by a constant unknown mean assumption in Section \ref{subsec:general}. 
For $a=1,\cdots, p$, we let $z_a$ have stationary covariance function $K_a:\R^d\to\R$ with $\Cov(z_a(\s), z_a(\s+\h))=K_a(\h).$

\begin{condition}\label{c:4}
A fixed $\delta>0$ exists such that, for all $n\in\mathbb{N}$ and, for all $i\neq j, i,j=1,\cdots,n, ||\s_i-\s_j||\ge\delta.$
\end{condition}

\begin{condition}\label{c:5}
Fixed $\beta>0$ and $\alpha>0$ exist such that, for all $\y\in\R^d$ and, for all $a=1,\cdots,p$
\[
|K_a(\y)|\le\frac{\beta}{1+||\y||^{d+\alpha}}.\]
\end{condition}

\begin{condition}\label{c:6}
Assuming Condition \ref{c:5} holds, then for the same $\beta>0$ and $\alpha>0$, we have, for $\ell=1,\cdots,k$
\[
|f_{\ell}(\y)|\le\frac{\beta}{1+||\y||^{d+\alpha}}.\]
\end{condition}

\begin{condition}\label{c:7}
For $\ell=1,\cdots,k$, we have
\[
\liminf_{n\to\infty}F_{n,f_{\ell}}>0.\]
\end{condition}

\begin{condition}\label{c:8}
For all $\ell, \ell^{'}=1,\cdots,k, \ell\neq \ell^{'}, f_{\ell}(\y)f_{\ell^{'}}(\y)=0$ for all $\y\in\R^d$.
\end{condition}

Condition \ref{c:4} implies that we are dealing with the increasing-domain asymptotic framework. For examples, see \cite{cressie2015statistics} for an introduction and \cite{bevilacqua2012estimating} for recent developments.
Condition \ref{c:5} holds for all the standard covariance functions in spatial statistics (\cite{cressie2015statistics}).
A typical example of a function $f\in\{f_1,\cdots,f_k\}$ for which Conditions \ref{c:2} and \ref{c:6} are satisfied is the ``ring" kernel:
\begin{equation}\label{eq:ring_kernel}
R(r_1, r_2)(\s)=I(r_1<||\s||\le r_2),
\end{equation}
with $0<r_1<r_2<\infty$. Condition \ref{c:7} is mild and simply requires that for $\ell=1,\cdots,k$, the number of pairs of observation locations $\s_i, \s_j, i, j=1,\cdots,n$, for which $f_{\ell}(\s_i-\s_j)$ is non-zero is not negligible when compared with $n$. Condition \ref{c:8} means that the supports of the kernels are disjoint. This enables us to have a simple and elegant asymptotic distribution of the test statistic.

Our first main result is on the asymptotic null distribution of our test statistic $t_q$.

\begin{proposition}\label{prop:asymt}
Assume that Conditions \ref{c:1}-\ref{c:8} hold. Then, as $n\to\infty$,
\[
t_q\xrightarrow{d}\chi^2_{k(p-q)(p-q+1)/2}.\]
\end{proposition}

In the next proposition, we show that when considering the same normalization as that considered by \cite{bachoc2018spatial} for the local covariance matrices, and when removing the assumption of disjoint kernel supports, we still obtain an asymptotic distribution of the test statistic as the distribution of the squared Euclidean norm of a Gaussian vector. In this proposition, we consider a metric $d_w$ generating the topology of weak convergence on the set of Borel probability measures on Euclidean spaces (e.g., \cite{dudley2018real}, p. 393).

\begin{proposition}\label{prop:dw}
Assume that Conditions \ref{c:1}-\ref{c:6} hold. Let the test statistic $\tilde{t}_q$ be defined as $t_q$, with $\wM(f)$ replaced by
\[
\widetilde{\M}(f)=\n\sum_{j=1}^nf(\s_i-\s_j)\x(\s_i)\x(\s_j)^T\]
for $f\in\{f_1,\cdots,f_k\}$. Let $\mathcal{L}_{\tilde{t}_q, n}$ be the distribution of the test statistic $\tilde{t}_q$, and let $\mathcal{L}_{\V,n}$ be the distribution of 
$
\sum_{\ell=1}^k\sum_{a,b=1}^{p-q}\V^2_{\ell,a,b},
$
where $(\V_{\ell,a,b})_{\ell=1,\cdots,k,a,b=1,\cdots,p-q}$ is a Gaussian vector with mean vector $\bm{0}$ and with covariance matrix defined by
\[
\Cov(\V_{\ell,a,b}, \V_{\ell', a', b'})=\frac{1}{2}F_{n, f_{\ell}, f_{\ell'}}(I(a=a')I(b=b')+I(a=b')I(b=a'))\]
with
\[
F_{n, f_{\ell}, f_{\ell'}}=\frac{1}{n}\sum_{i,j=1}^nf_{\ell}(\s_i-\s_j)f_{\ell'}(\s_i-\s_j),\]
for $\ell, \ell'=1,\cdots,k$ and $a, b, a', b'=1,\cdots,p-q$. Then, as $n\to\infty$,
\[
d_w(\mathcal{L}_{\tilde{t}_q, n}, \mathcal{L}_{\V,n})\to 0.\]
\end{proposition}

In the following corollary, we show that if the supports of the kernels are disjoint, the test statistic converges to a weighted chi-squared distribution. See e.g. \cite{Bodenham-Adams:2016} for a presentation of the approximation procedures for this distribution.

\begin{corollary}\label{coro:dw}
Consider the setting of Proposition \ref{prop:dw} and assume additionally that Condition \ref{c:8} holds. Then, the limiting distribution $\mathcal{L}_{\V,n}$ in Proposition \ref{prop:dw} is equal to the distribution of 
\[
\sum_{\ell=1}^kF_{n,f_{\ell}}\mathcal{X}_{\ell}^2\]
where $\mathcal{X}_1^2,\cdots,\mathcal{X}_k^2$ are independent and are chi-squared distributed with $(p-q)(p-q+1)/2$ degrees of freedom.
\end{corollary}

\subsection{Regular Domain as a Special Example}
When the data are observed in a regular-grid domain, i.e., $\mathcal{S}\subseteq\mathbb{Z}^d$, the kernel functions can be based on the natural notion of a spatial neighborhood on the grid, which simplifies our technique.

A location $\s_0=(s_1, \cdots, s_d)\in\mathbb{Z}^d$ has $2d$ one-way lag-$h$ neighbors, $(s_1\pm h, \cdots, s_d), (s_1, s_2\pm h, \cdots, s_d), \cdots, (s_1, \cdots, s_{d-1}, s_d\pm h)$. For example, if $d=2$ and $h=1$, the 4 one-way lag-$1$ neighbors of $\s_0$ are ``left" $(s_1-1, s_2)$, ``right" $(s_1+1, s_2)$, ``up" $(s_1, s_2+1)$ and ``down" $(s_1, s_2-1)$. Therefore, we could define the one-way lag-$1$ population and sample local covariance matrices as
\begin{small}
\begin{equation}
\M =
\frac{1}{\sqrt{n\sum_{i=1}^n|\mathcal{N}_{\s_i}|}}
\sum_{i=1}^n \sum_{\s_j\in \mathcal{N}_{\s_i}}
\mathbb{E}\left(\x(\s_i)\x(\s_j)^T\right) 
\quad \text{and} \quad 
\wM =
\frac{1}{\sqrt{n\sum_{i=1}^n|\mathcal{N}_{\s_i}|}}
\sum_{i=1}^n\sum_{\s_j\in \mathcal{N}_{\s_i}}
\x(\s_i)  \x(\s_j)^T
\end{equation}
\end{small}
where, for $\x \in \mathbb{Z}^d$, 
\[
\mathcal{N}_{\x}=\{\s \in \{ \s_1 , \ldots , \s_n \} ; |\x - \s | = 1 \},
\]
with $|\boldsymbol{u}| = |u_1| + \dots + |u_d|$ for $\boldsymbol{u} = (u_1,\ldots,u_d) \in \mathbb{R}^d$.
The matrices $\M$ and $\wM$ are of the form $\M(f)$ and $\wM(f)$ in \eqref{eq:plc} and \eqref{eq:slc} for $f(\s) = I( ||\s|| = 1 )$, $\s \in \mathbb{R}^d$.

Similarly, if $d=2$ a location $\s_0 = (s_1,s_2) \in \mathbb{Z}^2$ has $4$ two-way lag-$1$ neighbors that are of the form $(s_1 \pm 1 , s_2 \pm 1)$.
In general, for $m,h \in \mathbb{N}$, $1 \leq m \leq d$, the $m$-way lag-$h$ population and sample local covariance matrices can be defined as:
\begin{small}
\begin{equation}
\M=\frac{1}{\sqrt{n\sum_{i=1}^n|\mathcal{N}^m_{h, \s_i}|}}\sum_{i=1}^n\sum_{\s_j\in \mathcal{N}^m_{h, \s_i}}
\mathbb{E}\left(\x(\s_i)\x(\s_j)^T\right) 
\quad \text{and} \quad
\wM=\frac{1}{\sqrt{n\sum_{i=1}^n|\mathcal{N}^m_{h, \s_i}|}}\sum_{i=1}^n\sum_{\s_j\in \mathcal{N}^m_{h, \s_i}} \x(\s_i)\x(\s_j)^T,
\end{equation}
\end{small}
where, for $\x \in \mathbb{Z}^d$,
\[
\mathcal{N}^m_{h, \x}
=
\{
\s \in \{ \s_1 , \ldots , \s_n \};
 \s=\psi_J(\x, \zeta_J(\x) + h \boldsymbol{v}),
~ \text{for some} ~
  J \in \mathcal{A}_m, \boldsymbol{v} \in \{ - 1 , 1 \}^m
\},
\]
with $\mathcal{A}_m=\{J=(i_1,\cdots,i_m)\in\mathbb{N}^m; 1 \leq i_1<\cdots<i_m\le d\}$,
 that is
$|\mathcal{A}_m|=\binom{d}{m}$
and for $J =(i_1,\cdots,i_m) \in\mathcal{A}_m, \y=(y_1, \cdots, y_m)\in\mathbb{Z}^m, \zeta_J(\x)=(x_{i_1},\cdots,x_{i_m}), \psi_J(\x, \y)=(x_1, \cdots, x_{i_1-1}, y_1, x_{i_1+1}, \cdots, x_{i_m-1}, y_m, x_{i_m+1}, \cdots, x_d)$.

In general, in Equations~\eqref{eq:plc} and \eqref{eq:slc} the $m$-way lag-$h$ population and sample local covariance matrices $\M$ and $\wM$ can also be written in the form  $\M(f)$ and $\wM(f)$ , with $f(\s) = I( \s \in \{ -h,0,h \}^d, |\s| = hm )$, $\s \in \mathbb{R}^d$.

Consequently, for the similarly defined test statistic $t_q$, the same limiting conclusions can be derived. By exploiting the neighborhood structure in the regular domain case, we can also shorten the computation time for other techniques, such as the proposed asymptotic test or spatial bootstrap, with the greatest time improvement being achieved for the latter (see Sections \ref{subsec:bootstrap} and \ref{subsec:sim2}).

\subsection{Estimation of the Number of Signal Components}

In this section, we investigate an estimator of the signal number $q$ based on the asymptotic tests. Now we wish to test the null hypothesis, for $r\in\{0,\cdots,p-1\}$, 
\begin{center}
$H_{0r}:$ There are exactly $p-r$ white noise processes in $\z$.
\end{center}
This hypothesis states that the signal dimension is $r$. Similar to Section \ref{Sec:Setup}, for $r = 0,\ldots,p-1$,  we can partition, for $\ell = 1 , \ldots , k$,
\[
\widehat \D_{\ell} =
\left(
\begin{array}{cc}
\widehat \D_{\ell,rr} & \widehat \D_{\ell,r-r}\\
\widehat \D_{\ell,-rr} & \widehat \D_{\ell,-r-r}\\
\end{array}
\right)
\]
where the block $\widehat{\D}_{\ell,rr}$ has size $r \times r$ and the block $\widehat \D_{\ell,-r-r}$ has size $(p - r) \times (p-r)$. Then, consider then the test statistic:
\[
t_r = \frac{n}{2} \sum_{\ell=1}^k ||\widehat \D_{\ell,-r-r}||^2.
\]

We now use the test statistic $t_r, r=0,1,\cdots, p-1$ for the estimation problem and derive a number of useful limiting properties in the following proposition.

\begin{proposition} \label{prop:power}
Assume the same conditions as in Proposition \ref{prop:asymt}. Then, \begin{itemize}
\item If $r \geq q$, then $t_r$ is bounded in probability.
\item If $r < q$, then there exists a fixed $b >0$ such that $ t_r / n \geq b + o_p(1)$.
\end{itemize} 
\end{proposition}

A consistent estimate $\hat{q}$ of the unknown signal dimension $q\le p-1$ can then be based on the test statistic $t_r$ as follows.

\begin{proposition} \label{prop:estimating:q}
Assume the same conditions as in Proposition \ref{prop:asymt}.
Let $(c_n)_{n \in \mathbb{N}}$ be a sequence of positive numbers such that $c_n \to \infty$ and $c_n = o(n)$ as $n \to \infty$. Let
\[
\hat{q} = \min \left\{ r \in \{1,\dots,p-1\} \,|\,t_r \leq c_n \right\},
\]
 with the convention $\min \emptyset = p$. Then, $\hat{q} \to q$ in probability as $n \to \infty$.
\end{proposition}

Specifying the sequence $c_n$ is not obvious in practice.
Nevertheless, an estimator $\hat{q}$ can also be found by applying a suitable strategy to perform successive tests.
Later, in the simulations we will always test for simplicity at the same significance level and apply a divide-and-conquer strategy for the testing.

\subsection{General Mean}\label{subsec:general}

The previous results were derived under the assumption that $\mathbb{E}(\z(\s))=\0$. In the next proposition, we show that the conclusions of Propositions \ref{prop:asymt}, \ref{prop:dw}, \ref{prop:power}, and \ref{prop:estimating:q} and of Corollary \ref{coro:dw} are unchanged when $\z$ has a non-zero unknown constant mean function and when the observations are empirically centered for the computation of the local covariance matrices.

\begin{proposition}\label{prop:nonzerom}
Assume that for $a=1,\cdots,p, z_a$ has constant mean function $\mu_a\in\R$. Let, for $f\in\{f_1,\cdots,f_k\},$
\begin{equation}\label{eq:Mbar}
\widebar{\M}(f)=\frac{1}{n\sqrt{F_{n,f}}}\sum_{i=1}^n\sum_{j=1}^nf(\s_i-\s_j)(\x(\s_i)-\bar{\x})(\x(\s_j)-\bar{\x})^T \quad \text{and} \quad \widebar{\M}(f_0)=\n(\x(\s_i)-\bar{\x})(\x(\s_i)-\bar{\x})^T,
\end{equation}
with $\bar{\x}=(1/n)\sum_{i=1}^n\x(\s_i)$.\\
Then, the conclusions of Propositions \ref{prop:asymt}, \ref{prop:dw}, \ref{prop:power}, and \ref{prop:estimating:q} as well as that of Corollary \ref{coro:dw} still hold under the same assumptions, except that $\wM(f)$ is everywhere replaced by $\widebar{\M}(f)$.
\end{proposition}

For the remainder of the paper, we assume that the mean is unknown.

\subsection{Bootstrap Tests for Dimension}\label{subsec:bootstrap}

The above derived noise dimension test based on the large sample behavior of the introduced test statistic is efficient to compute, but a large sample size may be needed for the finite sample level to match the asymptotic one. As an alternative for smaller sample sizes, we can formulate noise dimension tests based on the bootstrap. These latter tests also have the benefit of not requiring $\z$ to be Gaussian.

In its original form, the bootstrap is a non-parametric tool for estimating the distribution of an estimator or test statistic by re-sampling from the empirical cumulative distribution function (ecdf) of the sample at hand. It has had good performance in many statistical problems by theoretical analysis as well as simulation studies and applications to real data. See \cite{chernick2011bootstrap} or \cite{lahiri} for a more detailed discussion.

Again, we assume that the observed random field is following the SBSS model given by Equation~\eqref{eq:model} and want to test $H_{0r}$ given an SBSS solution of Equation~\eqref{eq:hg} for a certain kernel setting and the corresponding test statistic seen in Equation~\eqref{eq:t}. In the following, we formulate a method for re-sampling from the distribution of Model \eqref{eq:model} by respecting the null hypothesis $H_{0r}$. In line with the ideas presented by \cite{bootstrap_ts} this is achieved by leaving the hypothetical signal part of the estimated latent field $\hat{\z}(\bo s) = \hat{ \G} \x(\bo s) $ untouched and manipulating only the hypothetical noise parts $(\hat{\z}(\bo s))_i$ for $i = r+1,\dots,p$ and all $\bo s \in \{\s_1 , \ldots , \s_n  \}$ in one of the following ways.

\textbf{Parametric:} Here, it is assumed that each noise part is independent and identically distributed (iid) Gaussian, as is usual for white noise processes. This leads to bootstrap samples $(\z^* (\bo s) )_{i} \sim N(0,1)$ for $i = r+1,\dots,p$ and corresponding to each $\bo s \in \{\s_1 , \ldots , \s_n  \}$.

\textbf{Permute:} Here, we assume that each noise component is still iid but that it does not necessarily follow a Gaussian distribution. Therefore, bootstrap samples are drawn from the ecdf of the joint noise components: $(\z^* (\bo s) )_{i} \sim \text{ecdf}((\hat{\z}(\bo s_1)^\top)_{\hat{w}}, \dots, (\hat{\z}(\bo s_n)^\top)_{\hat{w}})$, with $i=r+1,\ldots,p$, $\bo s \in \{ \s_1 , \ldots , \s_n \}$ and where $\hat{w}$ denotes the noise components ($r+1$ to $p$) of $\hat{\z}$.

\begin{algorithm}
	Set the number of resamples $B$, the observed sample $\X = (\x(\s_1), \ldots, \x(\s_n))^\top$, the flag $spatial \_ resampling$ and optionally the block size $m$;\;\\
	Compute the SBSS solution and get $\hat{\bo \Gamma}$ and $\hat{\bo Z} = (\hat{ \G} \X^\top)^\top $ and compute test statistic $t=t_r(\X)$;\;\\
	\For{$k \in \{ 1, \ldots , B \}$}{
	 Replace the last $p-r$ columns of $\hat{\bo Z}$ by either a parametric or bootstrap sample to get $\bo Z^{*k}$;\;\\
	\If{spatial\_resampling = TRUE}{Replace $\bo Z^{*k}$ by a full spatial bootstrap sample. See text for details.;}
	Compute $\bo X^{*k} \leftarrow \hat{\bo \Gamma} \bo Z^{*k} $ and $t^k \leftarrow t_r(\bo X^{*k})$;\;
}
	Return the $p$-value: $[\#(t^k \geq t) + 1]/(B + 1)$;\;
	\caption{Testing $H_{0r} : q=r$}
	\label{alg::boot}
\end{algorithm}

After replacing the hypothetical noise part by a bootstrap sample in one of the former ways, the goal of sampling from Model \eqref{eq:model} under $H_{0r}$ is achieved. However, so far the uncertainty of estimating the signal has not been considered in the bootstrap test. Therefore, an optional second step in the whole re-sampling procedure is devoted to drawing a spatial bootstrap sample from the already manipulated sample as follows. We suggest the application of spatial bootstrapping as discussed in \cite{lahiri}, and in the following we summarize the main ideas. Let us recall that the set of sampling sites $\mathcal{C} = \{ \s_1 , \ldots , \s_n \}$ lies inside the $d$-dimensional spatial domain $\mathcal{S}$, which can be viewed as the ``sample region'' and hence $\mathcal{C} \subseteq \mathcal{S} \subseteq \mathbb{R} ^ d$. $\mathcal{S}$ is divided into non-overlapping blocks of size $m^d$ that lie partially in $\mathcal{S}$, formally $\mathcal{B} = \{ b_{\bo i} =  (\bo i + (0, 1] ^ d)m \cap \mathcal{S} : (\bo i + (0, 1] ^ d)m \cap \mathcal{S} \neq \emptyset, \bo i \in \mathbb{Z}^d \}$, and overlapping blocks that lie fully in $\mathcal{S}$, written as $\mathcal{B}_{bs} = \{ b_{\bo j} = \bo j + (0, 1] ^ d m : \bo j + (0, 1] ^ d m \subseteq \mathcal{S}, \bo j \in \mathbb{Z}^d \}$. The bootstrapped spatial domain $\mathcal{S}^*$ is formed by replacing each block $b_{\bo i} \in \mathcal{B}$ with a randomly with replacement sampled block $b_{\bo j} \in \mathcal{B}_{bs}$ that is trimmed to the shape of $b_{\bo i}$ by $b_{\bo j} \cap (b_{\bo i} - \bo i m + \bo j)$. Hence, the trimmed version of $b_{\bo j}$ remains at the original location of $b_{\bo j}$, while the shape changes to that of $b_{\bo i}$, taking care of the boundary blocks that do not fully lie within $\mathcal{S}$. Finally, the bootstrapped version of the random field writes as $\bo z^* = \{\bo z(\bo s) : \bo s \in \mathcal{S}^* \cap \mathcal{C} \}$. Note that in each spatial bootstrap iteration, the shape of $\mathcal{S}^*$ and therefore the bootstrapped sampling sites differ. This in turn makes the computation of the local covariance matrices a demanding task, as it relies on the distances between all sampling sites, which need to be newly computed in each iteration. For regular data, this can be avoided by using a slightly different bootstrap regime as follows.

\cite{lahiri_regular} have suggested a slightly different approach for sampling sites located on a regular grid, meaning that the sampling sites satisfy $\{ \bo s_1, \dots, \bo s_n \} \subseteq \mathcal{S} \cap \mathbb{Z}^d$. Again, the domain $\mathcal{S}$ is divided into blocks of size $m ^ d$ that are either non-overlapping or overlapping but lie completely inside $\mathcal{S}$, leading to $\mathcal{B} = \{ (\bo i + (0, 1] ^ d)m : (\bo i + (0, 1] ^ d)m \subseteq \mathcal{S}, \bo i \in \mathbb{Z}^d \}$ and $\mathcal{B}_{bs}$, as defined above. The key difference is that the bootstrap sample is drawn at the level of the random field values, whereas the former bootstrap version operates at the level of the spatial domain. Specifically, for each block $b_{\bo i} \in \mathcal{B}$ the values $\{\z(\bo s) : \bo s \in b_{\bo i} \cap \mathbb{Z}^d\}$ are replaced by $\{\z(\bo s) : \bo s \in b_{\bo j} \cap \mathbb{Z}^d \}$ for a randomly with replacement chosen block $b_{\bo j} \in \mathcal{B}_{bs}$. This procedure keeps the bootstrapped spatial domain and sampling sites equal in all iterations, namely the unison of all blocks from $\mathcal{B}$. This in turn simplifies the computation of local covariance matrices, as only the random field values change. We will compare the computation times of the former two approaches in the simulation study presented in Section~\ref{subsec:sim2}.

Algorithm~\ref{alg::boot} summarizes the formerly discussed bootstrap strategy to test for one specific value of signal dimension $r$. To estimate the signal dimension, a sequence of tests for different signal dimensions $r$ at a given significance level $\alpha$ are carried out. A number of different test sequences are possible, but we rely on a divide-and-conquer strategy outlined in Algorithm~\ref{alg::div_con}. Here, the $test \_ function$ could be either one of the bootstrap test variants seen in Algorithm~\ref{alg::boot} or the asymptotic test outlined above.


\begin{algorithm}
	Set $lower$, $upper$ and $\alpha$;\;\\
	$middle = \lfloor (upper-lower)/2 \rfloor$;\;\\
	\While{$(middle != lower) ~ \& \& ~ (middle != upper)$}{
	 $p = test \_ function(r=middle)$;\;\\
	\If{$p < alpha$}{$lower = middle$;}
	\Else{$upper = middle$;}
	$middle = \lfloor (upper-lower)/2 \rfloor$;\;	}
	Return $\hat{q} = middle + 1$;\;
	\caption{Divide and Conquer}
	\label{alg::div_con}
\end{algorithm}

\section{Simulation}\label{Sec:Sim}

To validate the performance of the methods introduced above, we carried out three extensive simulation studies in R 3.6.1 (\cite{r_language}) with the help of the packages SpatialBSS (\cite{spatialbss_package}), JADE (\cite{jade_package}), sp (\cite{sp_package}), raster (\cite{raster_package}), gstat (\cite{gstat_package}) and RandomFields (\cite{randomfields_package}).

\subsection{Simulation Study 1: Hypothesis Testing}\label{sec::sim_hypothesis}
\begin{figure}
\centering
    \begin{minipage}[t]{0.32\linewidth}
        \centering
        \includegraphics[scale=0.5]{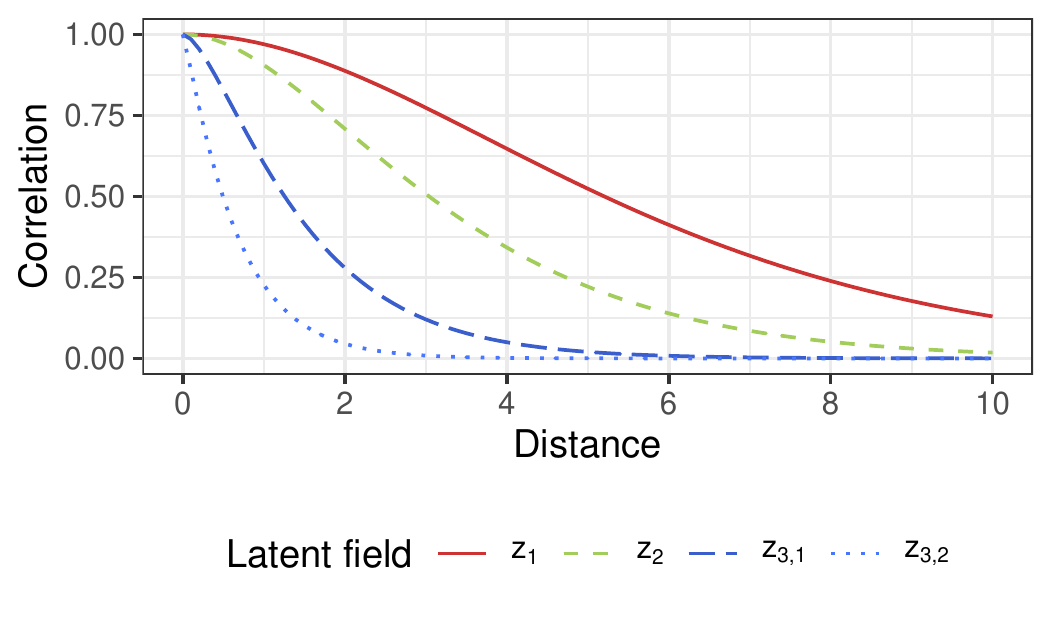}
    \end{minipage}
    \hfill
    \begin{minipage}[t]{0.32\linewidth}
        \centering
        \includegraphics[scale=0.5]{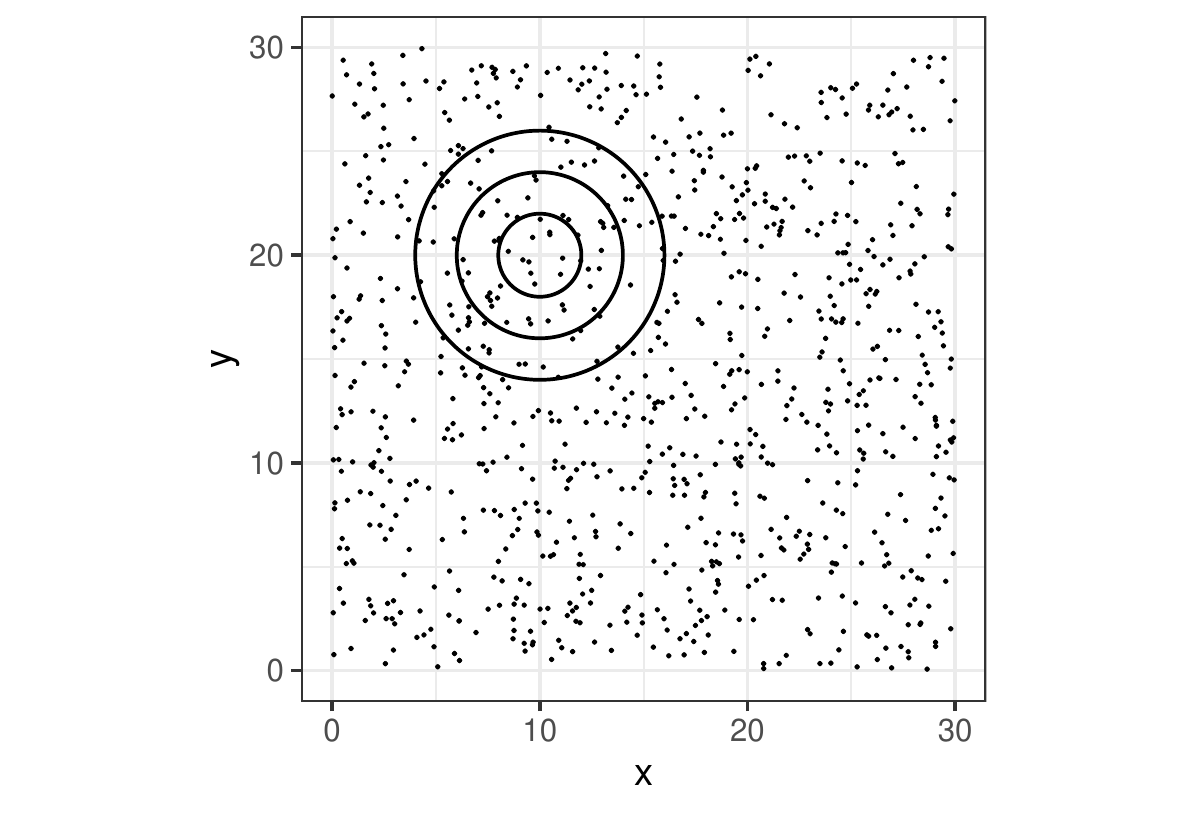}
    \end{minipage}
    \hfill
    \begin{minipage}[t]{0.32\linewidth}    
        \centering
        \includegraphics[scale=0.5]{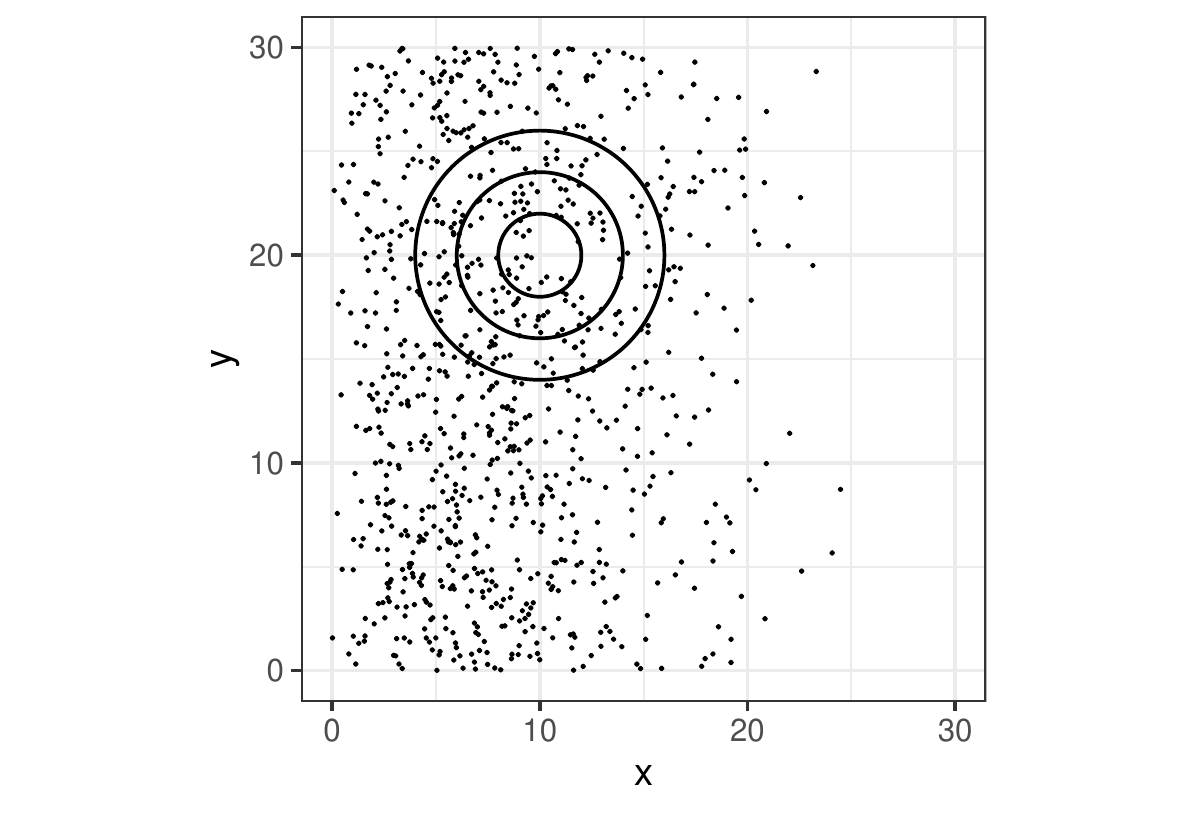}
    \end{minipage}
    \caption{
Left: Mat\'ern correlation functions for model setting 1, which consists of the signal random field $(z_1,z_2,z_{3,1})$ with parameters $(\nu, \phi) \in \{ (3,2), (2,1.5), (1,1) \}$ and model setting 2 formed by the signal random field $(z_1,z_2,z_{3,2})$ with parameters $(\nu, \phi) \in \{(3,2), (2,1.5), (0.6,0.6)\}$.
Middle and right: uniform (middle) and skewed (right) coordinate sample pattern for a spatial domain of size $30 \times 30$ with three circles of radii $(2,4,6)$ representing ring kernel functions.}
        \label{fig:coordinates}
\end{figure}

In this part of the simulation, we explored the performance of hypothesis testing. For all the following simulations, we considered the SBSS model, as shown in Equation~\eqref{eq:model}, where without loss of generality we set $\mu_a = 0$ for $a=1,\dots,p$ and assume the mean to be unknown. For the latent signal part we used two different three-variate random field model settings. Therefore, the true dimension is always $q=3$. All the random fields followed a Mat\'ern correlation structure, and the $a$-th random field $z_a$ thus had its covariance function value at $\boldsymbol{u}, \boldsymbol{v} \in \mathbb{R}^d$, given by:
\[
K_a(h; \nu, \phi) = \frac{1}{2 ^ {\nu - 1} \Gamma (\nu)} \left( \frac{h}{\phi} \right) ^ \nu  K_\nu \left( \frac{h}{\phi} \right),
~ ~
h = || \boldsymbol{u} - \boldsymbol{v} ||,
\]
where $\nu >0$ is the shape parameter, $\phi >0$ is the range parameter, and $K_\nu$ is the modified Bessel function of second kind with shape parameter $\nu$. The parameters used were $(\nu, \phi) \in \{(3,2), (2,1.5), (1,1) \}$ and $\{(3,2), (2,1.5), (0.6,0.6) \}$ for model setting 1 and 2 respectively, which are depicted in Figure~\ref{fig:coordinates}. Model setting 2 can be viewed as a low-dependence version of model setting 1. The noise part always consists of iid samples drawn from $N_{2}(\bo 0, \I_{2})$, leading to a total latent field dimension of $p=5$ for both model settings. As SBSS is affine equivariant (for details see \cite{bachoc2018spatial} and the appendix) we chose the mixing matrix to be the identity matrix, i.e., $\Om = \bf I_5$, without loss of generality.

We focused on the squared spatial domains $[0,n] \times [0,n]$ (also written in the following as $n \times n$) of different sizes $n \in \{30, 40, 50, 60\}$. For a given domain, we considered two different sample location patterns: uniform and skewed. For the uniform pattern, $n ^ 2$ pairs of $(x,y)$-coordinates were randomly drawn from a uniform distribution $U(0,1)$ and then multiplied by $n$, leading to a constant sampling location density over the entire domain. We followed the same approach for the skewed pattern, with the only difference being that the $x$ coordinate values were drawn from a beta distribution $\beta (2,5)$, resulting in a denser arrangement of samples in the left half of the domain.

For the local covariance matrices \eqref{eq:slc}, we used two different kernel function settings. Kernel setting 1 used only one ring kernel function \eqref{eq:ring_kernel} with parameters $(r_1,r_2) = (0,2)$, while kernel setting 2 used three ring kernel functions with parameters $(r_1,r_2) \in \{(0, 2), (2, 4), (4, 6)\}$. Figure~\ref{fig:coordinates} depicts a simulation example for each of the uniform and skewed coordinate patterns, where the circles represent the different ring kernel radii.

For each of the four simulation settings, we carried out 2000 repetitions, and in each repetition we tested three different null hypothesis ($H_{02}$, $H_{03}$, and $H_{04}$) with the following five test approaches: asymptotic test (Asym), noise bootstrapping with option parametric (Param), noise bootstrapping with option permute (Perm), full spatial bootstrapping with option parametric (Sp Param), and full spatial bootstrapping with option permute (Sp Perm). For all bootstrap approaches, we fixed the number of re-samples to be $B=200$, and for the full spatial bootstrap, the block size was equal to $m=10$. 

Rejection rates based on a significance level of $\alpha = 0.05$ for all simulation settings are presented in Tables~\ref{tab:model_1} and \ref{tab:model_2}. Overall, all the test methods appeared to maintain the expected rejection rates, which were $1.00$ for $H_{02}$, $0.05$ for $H_{03}$, and $<0.05$ for $H_{04}$ based on $\alpha = 0.05$. Only for small samples sizes ($30 \times 30$) did the asymptotic test show a too small rejection rate for kernel setting 2 and the skewed sample location pattern. Thus, for practical applications, smaller numbers of kernel functions might be preferable for the asymptotic test. For bootstrapping, the full spatial variants and those relying only on manipulating the hypothetical noise part performed equally well. Considering the computation time, the latter bootstrap variant might be preferable, as explored in more detail in Section \ref{subsec:sim2}.

\begin{table}[h]
\caption{Rejection rates for model setting 1  at level $\alpha=0.05$.}
\label{tab:model_1}
\centering
\scalebox{0.8}{
\begin{tabular}{clcccccccccccc}
  \toprule
  &  &  \multicolumn{6}{c}{Uniform} & \multicolumn{6}{c}{Skew} \\
   \cmidrule(r){3-8} \cmidrule(l){9-14}
   &  & \multicolumn{3}{c}{Kernel Setting 1} & \multicolumn{3}{c}{Kernel Setting 2} & \multicolumn{3}{c}{Kernel Setting 1} & \multicolumn{3}{c}{Kernel Setting 2}\\
   \cmidrule(r){3-5} \cmidrule(l){6-8} \cmidrule(l){9-11} \cmidrule(l){12-14}
Domain & Method & $H_{02}$ & $H_{03}$ & $H_{04}$ & $H_{02}$ & $H_{03}$  & $H_{04}$ & $H_{02}$ & $H_{03}$ & $H_{04}$ & $H_{02}$ & $H_{03}$  & $H_{04}$ \\ 
  \midrule
\multirow{5}{*}{$30 \times 30$} & Asym & 1.000 & 0.041 & 0.006 & 1.000 & 0.042 & 0.007  & 1.000 & 0.042 & 0.004 & 1.000 & 0.029 & 0.003\\ 
   & Sp Param & 1.000 & 0.048 & 0.006 & 1.000 & 0.058 & 0.001 & 1.000 & 0.059 & 0.004 & 1.000 & 0.051 & 0.000\\ 
   & Sp Perm & 1.000 & 0.050 & 0.006 & 1.000 & 0.059 & 0.000 & 1.000 & 0.058 & 0.004 & 1.000 & 0.052 & 0.001\\ 
   & Param & 1.000 & 0.042 & 0.006 & 1.000 & 0.044 & 0.006 & 1.000 & 0.050 & 0.008 & 1.000 & 0.039 & 0.005\\ 
   & Perm & 1.000 & 0.045 & 0.008 & 1.000 & 0.051 & 0.006 & 1.000 & 0.049 & 0.008 & 1.000 & 0.035 & 0.005\\ \midrule
  \multirow{5}{*}{$40 \times 40$} & Asym & 1.000 & 0.055 & 0.004 & 1.000 & 0.048 & 0.005 & 1.000 & 0.045 & 0.002 & 1.000 & 0.040 & 0.005\\ 
   & Sp Param & 1.000 & 0.056 & 0.005 & 1.000 & 0.066 & 0.000 & 1.000 & 0.056 & 0.003 & 1.000 & 0.064 & 0.002\\ 
   & Sp Perm & 1.000 & 0.063 & 0.005 & 1.000 & 0.061 & 0.000 & 1.000 & 0.055 & 0.004 & 1.000 & 0.065 & 0.002 \\ 
   & Param & 1.000 & 0.052 & 0.007 & 1.000 & 0.055 & 0.003 & 1.000 & 0.050 & 0.007 & 1.000 & 0.048 & 0.005\\ 
   & Perm & 1.000 & 0.056 & 0.007 & 1.000 & 0.052 & 0.004 & 1.000 & 0.048 & 0.008 & 1.000 & 0.050 & 0.004\\ \midrule
  \multirow{5}{*}{$50 \times 50$} & Asym & 1.000 & 0.049 & 0.005 & 1.000 & 0.040 & 0.010 & 1.000 & 0.040 & 0.006 & 1.000 & 0.044 & 0.009\\ 
   & Sp Param & 1.000 & 0.052 & 0.004 & 1.000 & 0.053 & 0.002 & 1.000 & 0.047 & 0.006 & 1.000 & 0.064 & 0.002\\ 
   & Sp Perm & 1.000 & 0.050 & 0.005 & 1.000 & 0.053 & 0.002 & 1.000 & 0.045 & 0.005 & 1.000 & 0.061 & 0.002\\ 
   & Param & 1.000 & 0.052 & 0.007 & 1.000 & 0.049 & 0.007 & 1.000 & 0.042 & 0.007 & 1.000 & 0.054 & 0.007\\ 
   & Perm & 1.000 & 0.050 & 0.008 & 1.000 & 0.050 & 0.006  & 1.000 & 0.042 & 0.010 & 1.000 & 0.054 & 0.008\\ \midrule
  \multirow{5}{*}{$60 \times 60$} & Asym & 1.000 & 0.052 & 0.006 & 1.000 & 0.048 & 0.010 & 1.000 & 0.044 & 0.004 & 1.000 & 0.045 & 0.004\\ 
   & Sp Param & 1.000 & 0.056 & 0.006 & 1.000 & 0.058 & 0.003 & 1.000 & 0.048 & 0.005 & 1.000 & 0.060 & 0.000\\ 
   & Sp Perm & 1.000 & 0.055 & 0.007 & 1.000 & 0.057 & 0.002 & 1.000 & 0.052 & 0.004 & 1.000 & 0.058 & 0.000\\ 
   & Param & 1.000 & 0.049 & 0.009 & 1.000 & 0.054 & 0.006 & 1.000 & 0.043 & 0.006 & 1.000 & 0.048 & 0.004\\ 
   & Perm & 1.000 & 0.053 & 0.009 & 1.000 & 0.050 & 0.008 & 1.000 & 0.046 & 0.006 & 1.000 & 0.048 & 0.004\\ 
  \bottomrule
\end{tabular}
}
\end{table}

\begin{table}[h]
\caption{Rejection rates for model setting 2 at level $\alpha=0.05$.}
\centering
\label{tab:model_2}
\scalebox{0.8}{
\begin{tabular}{clcccccccccccc}
  \toprule
    &  &  \multicolumn{6}{c}{Uniform} & \multicolumn{6}{c}{Skew} \\
     \cmidrule(r){3-8} \cmidrule(l){9-14}
   &  & \multicolumn{3}{c}{Kernel Setting 1} & \multicolumn{3}{c}{Kernel Setting 2} & \multicolumn{3}{c}{Kernel Setting 1} & \multicolumn{3}{c}{Kernel Setting 2}\\
   \cmidrule(r){3-5} \cmidrule(l){6-8} \cmidrule(l){9-11} \cmidrule(l){12-14}
Domain & Method & $H_{02}$ & $H_{03}$ & $H_{04}$ & $H_{02}$ & $H_{03}$  & $H_{04}$ & $H_{02}$ & $H_{03}$ & $H_{04}$ & $H_{02}$ & $H_{03}$  & $H_{04}$ \\ 
  \midrule
\multirow{5}{*}{$30 \times 30$} & Asym & 1.000 & 0.051 & 0.005 & 1.000 & 0.052 & 0.004 & 1.000 & 0.048 & 0.006 & 1.000 & 0.033 & 0.003\\ 
   & Sp Param & 1.000 & 0.053 & 0.005 & 1.000 & 0.062 & 0.000 & 1.000 & 0.058 & 0.005 & 1.000 & 0.055 & 0.002\\ 
   & Sp Perm & 1.000 & 0.052 & 0.006 & 1.000 & 0.065 & 0.001 & 1.000 & 0.056 & 0.006 & 1.000 & 0.051 & 0.001\\ 
   & Param & 1.000 & 0.052 & 0.011 & 1.000 & 0.058 & 0.003 & 1.000 & 0.059 & 0.011 & 1.000 & 0.043 & 0.004\\ 
   & Perm & 1.000 & 0.048 & 0.011 & 1.000 & 0.060 & 0.002 & 1.000 & 0.061 & 0.012 & 1.000 & 0.044 & 0.003\\ \midrule
  \multirow{5}{*}{$40 \times 40$} & Asym & 1.000 & 0.060 & 0.004 & 1.000 & 0.052 & 0.005 & 1.000 & 0.050 & 0.004 & 1.000 & 0.038 & 0.007\\ 
   & Sp Param & 1.000 & 0.063 & 0.002 & 1.000 & 0.060 & 0.000 & 1.000 & 0.060 & 0.004 & 1.000 & 0.054 & 0.002\\ 
   & Sp Perm & 1.000 & 0.055 & 0.002 & 1.000 & 0.062 & 0.000 & 1.000 & 0.058 & 0.002 & 1.000 & 0.057 & 0.002\\ 
   & Param & 1.000 & 0.056 & 0.006 & 1.000 & 0.056 & 0.004 & 1.000 & 0.052 & 0.008 & 1.000 & 0.045 & 0.005\\ 
   & Perm & 1.000 & 0.058 & 0.005 & 1.000 & 0.053 & 0.004 & 1.000 & 0.054 & 0.006 & 1.000 & 0.045 & 0.005\\ \midrule
  \multirow{5}{*}{$50 \times 50$} & Asym & 1.000 & 0.045 & 0.004 & 1.000 & 0.047 & 0.004 & 1.000 & 0.044 & 0.005 & 1.000 & 0.044 & 0.004\\ 
   & Sp Param & 1.000 & 0.048 & 0.002 & 1.000 & 0.056 & 0.000 & 1.000 & 0.053 & 0.002 & 1.000 & 0.058 & 0.001\\ 
   & Sp Perm & 1.000 & 0.049 & 0.002 & 1.000 & 0.053 & 0.001 & 1.000 & 0.050 & 0.005 & 1.000 & 0.055 & 0.001\\ 
   & Param & 1.000 & 0.045 & 0.004 & 1.000 & 0.050 & 0.002 & 1.000 & 0.048 & 0.007 & 1.000 & 0.051 & 0.004\\ 
   & Perm & 1.000 & 0.044 & 0.007 & 1.000 & 0.048 & 0.003 & 1.000 & 0.046 & 0.009 & 1.000 & 0.052 & 0.004\\ \midrule
  \multirow{5}{*}{$60 \times 60$} & Asym & 1.000 & 0.048 & 0.004 & 1.000 & 0.059 & 0.008 & 1.000 & 0.047 & 0.004 & 1.000 & 0.042 & 0.006\\ 
   & Sp Param & 1.000 & 0.052 & 0.005 & 1.000 & 0.072 & 0.002 & 1.000 & 0.050 & 0.004 & 1.000 & 0.059 & 0.000\\ 
   & Sp Perm & 1.000 & 0.056 & 0.003 & 1.000 & 0.068 & 0.002 & 1.000 & 0.050 & 0.004 & 1.000 & 0.057 & 0.000\\ 
   & Param & 1.000 & 0.047 & 0.009 & 1.000 & 0.063 & 0.004 & 1.000 & 0.046 & 0.005 & 1.000 & 0.052 & 0.003\\ 
   & Perm & 1.000 & 0.048 & 0.010 & 1.000 & 0.063 & 0.006 & 1.000 & 0.048 & 0.005 & 1.000 & 0.050 & 0.005\\ 
   \bottomrule
\end{tabular}
}
\end{table}

\subsection{Simulation Study 2: Computation Time Comparison}\label{subsec:sim2}

In this simulation, we investigated the computation times for the various test methods. As an illustrative example, we again considered a five-variate latent random field with model setting 1 and bivariate Gaussian noise components. In addition, we kept the same spatial domain sizes, though the sampling sites were changed to be regular defined as $[0,n] \times [0,n] \cap \mathbb{Z}^2$. $H_{03}$ was tested using the five former mentioned test methods with the same number of bootstrap samples and block sizes. The key difference is that each test was carried out with code designed for irregular sample locations as well as code that takes into account simplifications made possible by the fact that the sample locations were regular (e.g., the simplified spatial bootstrap algorithm). Two ring kernel functions with parameters $(r_1, r_2) \in \{(0, 1), (1, \sqrt{2}) \}$ were considered for the irregular code, and kernels of the form $f(\s)=I(||\s||=h)$ with $h \in \{1, \sqrt{2} \}$ were considered for the regular code (one-way and two-way lag-$1$ local covariance matrices). This choice ensured that the same neighbors were selected for both versions of the code and thus that the qualitative results of the tests were equal up to random effects of the bootstrap sampling procedures. 

\begin{figure}[th]
	\centering
  \includegraphics[scale=0.5]{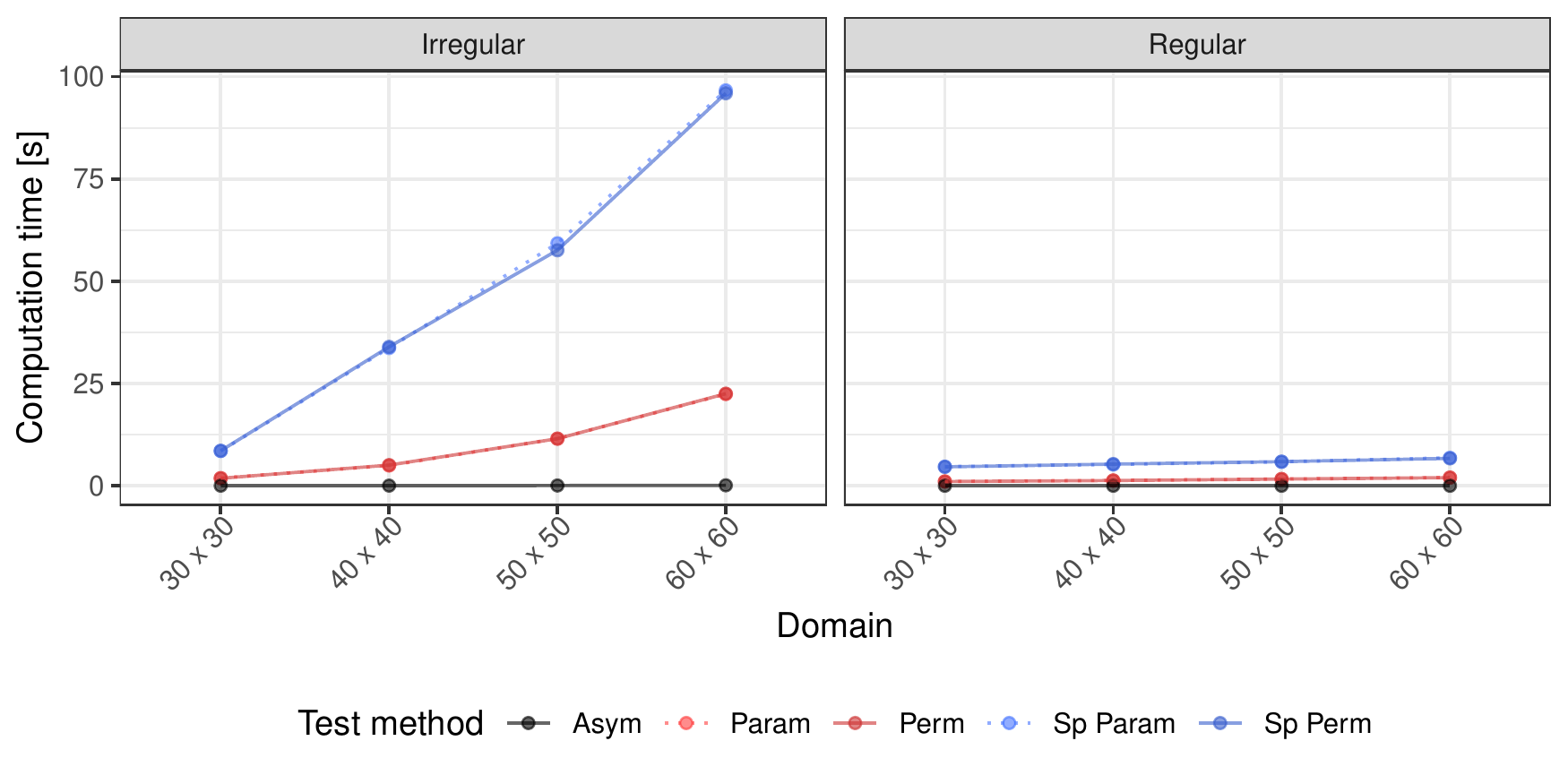}
  \caption{Median running times of the five different test methods for different domain sizes with regular sampling sites based on five simulation repetitions. Computations were carried out with code designed for regular and irregular sampling sites.}
  \label{fig::2k_running_times}
\end{figure}

Figure~\ref{fig::2k_running_times} shows the median computation time based on five simulation repetitions carried out on a Windows machine with an Intel i5 CPU. The computation times revealed that asymptotic tests are fastest, as the SBSS solution needs to be computed only once, whereas bootstrap algorithms compute the SBSS solution $B$ times.

Of greater interest is the overall difference in the computation time between regular and irregular code. This might be explained by the fact that the code for regular sampling sites does not rely on distances between sampling sites as the irregular code does. Specifically, the selection of neighbors for local covariance matrices can be implemented by shifting the coordinate system appropriately for the regular code, whereas in the irregular code this is based on looping over the distance matrix among all coordinates. This difference should also explain the different scaling of the computation time with increasing sample size, as looping through the distance matrix depends on the actual number of locations, while coordinate shifting does not.

Further, there was a larger computation time difference between the full spatial bootstrap and the one that manipulates only the hypothetical noise for the irregular code compared with the regular one. This might be the impact of the simplified spatial bootstrap variant for regular sampling sites. As explained above, for the irregular code the distance matrix has to be computed again for every new iteration because the spatial bootstrap changes sampling sites for each iteration, which is not the case for the regular code, for which the sampling sites remain equal for each bootstrap iteration.

Overall, this simulation strongly indicates that regular sampling sites should be computationally treated as such. In addition, considering the overall similar performances of the tests in the former simulation, the spatial bootstrapping step for the irregular data might be discarded, as it significantly increases the computation time.

\subsection{Simulation Study 3: Estimation of the Signal Dimension}\label{subsec:sim3}

The former simulations investigated only hypothesis tests for one specific value of the hypothetical signal dimension. In this section, we explore the use of hypothesis tests for signal dimension estimation. We considered the exact same simulation settings as in Section~\ref{sec::sim_hypothesis} but increased the dimension of the noise part to seven leading to a total latent random field dimension of $p=10$, while the true signal dimension remained $q=3$. Estimation of the signal dimension was based on the divide-and-conquer strategy described above. All hypothesis tests were carried out using the asymptotic test method and the parametric bootstrap without full spatial bootstrapping. This choice is justified by the similar performance in signal dimension testing of all bootstrap test variants and the fact that the full spatial bootstrap is computationally unfeasible for such a large simulation.

\begin{figure}[ht]
	\centering
  \includegraphics[scale=0.5]{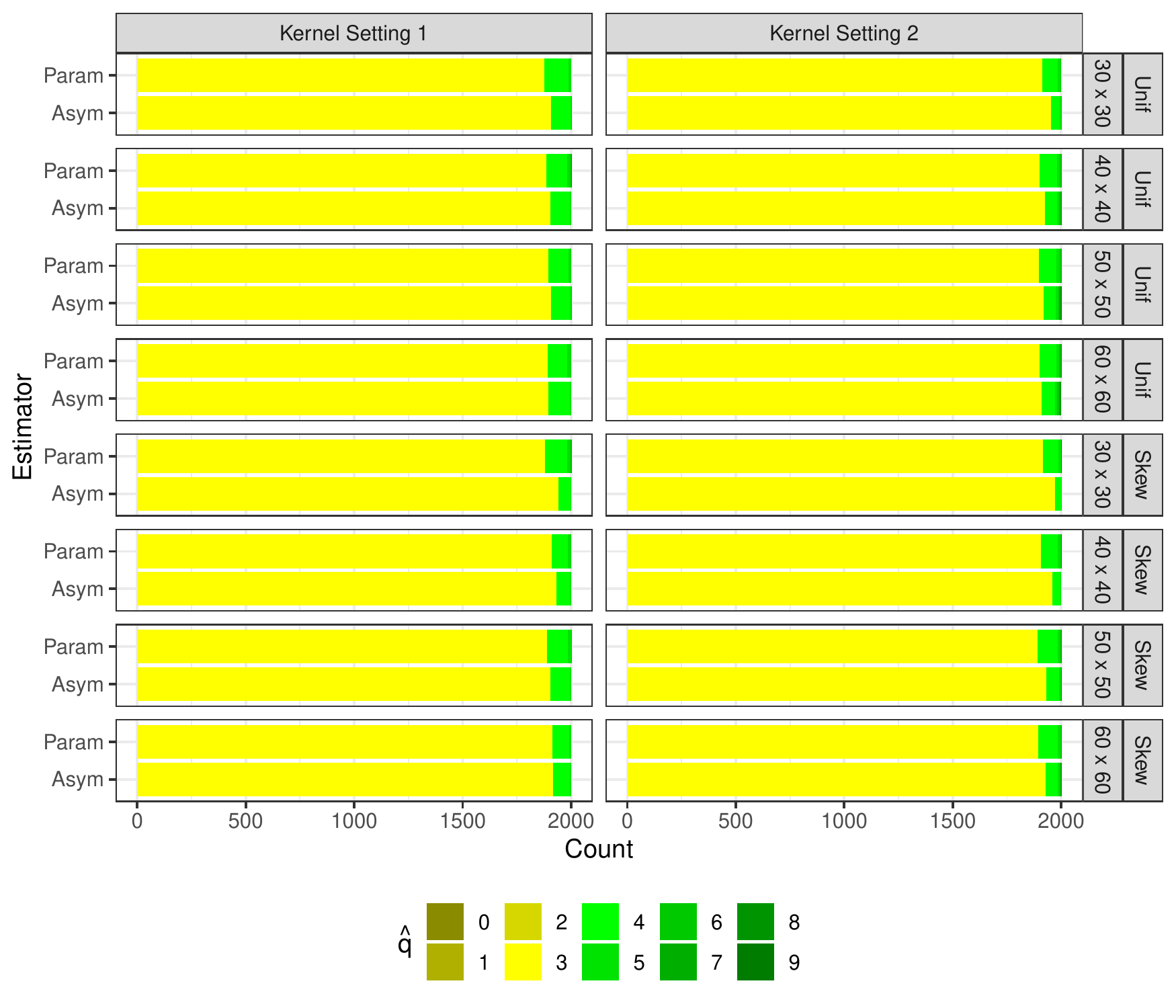}
  \caption{Frequencies of the estimated signal dimension for model setting 1.}
  \label{fig::est_1}
\end{figure}

\begin{figure}[ht]
	\centering
  \includegraphics[scale=0.5]{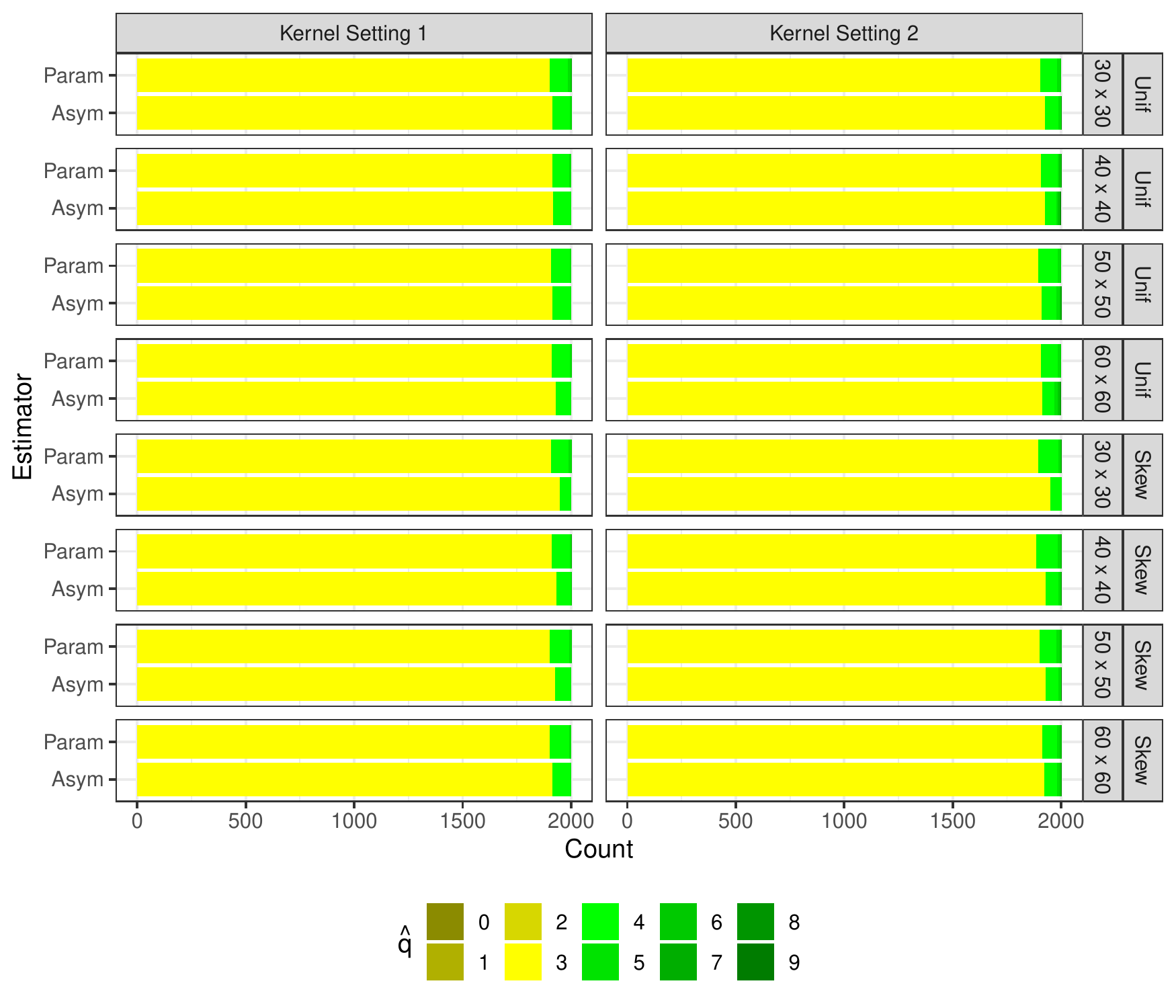}
  \caption{Frequencies of the estimated signal dimension for model setting 2.}
  \label{fig::est_2}
\end{figure}

Figures~\ref{fig::est_1} and \ref{fig::est_2} depict the estimated dimensions for 2000 simulation repetitions for a significance level of $\alpha=0.05$.
Overall, the estimation was highly accurate, with the estimated dimension being equal to the true one in approximately $95 \%$ of the cases.
 Interestingly, the signal dimension was never underestimated, while it was overestimated in approximately 100 of the simulation iterations, reflecting the significance level $\alpha=0.05$. For all settings, the asymptotic test performed better than the bootstrap test. This was especially true for low sample sizes, which is a counterintuitive result. However, it may be due to the fact that, as the former simulations show, for low sample sizes the asymptotic test never met the theoretical rejection rate, which is simply the significance level when the null is actually true for small sample sizes (Tables~\ref{tab:model_1} and \ref{tab:model_2}). Therefore, the true null is more often accepted leading to a better performance when estimating the signal dimension. 

\section{Data Example}\label{Sec:Ex}

In this section, we illustrate the application of our tests to a geochemical dataset. Specifically, we analyze samples of the amount of 31 elements in terrestrial moss, which were collected during the Kola project at 594 sites located in Norway, Finland, and Russia alongside the coast of the Barents Sea. The data is freely available in the R package \cite{statda_package} and is described in more detail by \cite{moss_book}. \cite{nordhausen2015blind} and \cite{bachoc2018spatial} have already considered the dataset in the context of SBSS where the goal of the former publication was to identify geochemical interpretative components. Indeed, six meaningful components were found based on an SBSS solution that used a single local covariance matrix with a ball kernel $B(\s) = I(||\s|| \leq 50 \text{km})$, chosen based on experts knowledge. The latter publication introduces a variety of kernel function options and expands SBSS for the use of more than one local covariance matrix. The authors also found that jointly diagonalizing more than one local covariance matrix yielded a more stable solution compared with using only one. For the moss dataset, an SBSS solution based on four local covariance matrices exclusively using ring kernel functions with parameters $(r_1,r_2) \in \{(0, 25), (25, 50), (50, 75), (75, 100)\}$ km was used, among others, where the first six latent field components showed high correlation with the one presented by \cite{nordhausen2015blind}. 

This was the starting point of our analysis, in which we estimated the signal dimension using the five test methods defined in Section \ref{sec::sim_hypothesis} based on the two kernel function settings described above. Note that, theoretically, the ball kernel function is not compatible with Condition~\ref{c:2}. Hence, we used exclusively ring kernel functions, leading to the two considered kernel settings: $(r_1,r_2)= (0, 25)$ (Setting 1) and $(r_1,r_2) \in \{(0, 25), (25, 50), (50, 75), (75, 100)\}$ (Setting 2) in km. Using the same methods as \cite{nordhausen2015blind} and \cite{bachoc2018spatial}, we pre-processed the data by performing an isometric log-ratio (ilr) transformation using pivot coordinates to respect the compositional nature of the geochemical data, which reduced the dimension of the data to $p=30$. The details of the ilr transformation and on compositional data in general can be found in \cite{comp_book}. For the bootstrap methods $B=200$ re-samples were drawn, and for the spatial bootstrap the domain was overlaid by a square grid with a side-length of 30 km. The block size was chosen to be $m=60$ km, meaning that every 30 km, a block of size 60 km was placed, forming the set of blocks from which the samples were drawn. All lengths were referred to the UTM zone 35N.

Table \ref{tab:kola} summarizes the estimated signal dimensions for the various test methods and kernel function settings based on a significance level of $\alpha = 0.05$. For both kernel settings, the dimension of the actual signal was approximately half of the original data dimension. This is a significant reduction in dimension, as further analysis (such as spatial prediction) needs to consider drastically fewer signal fields. Still, the estimated number of signals was higher than the number of meaningful geochemical components found by \cite{nordhausen2015blind}. The estimated number of signal dimensions for kernel setting 1 was always $\hat{q} = 14$, meaning it appears to be more stable (less sensitive to the choice of the test version) than estimations based on a greater number of ring kernel functions (setting 2), where the number ranged between $\hat{q} = 15$ and $18$. Thus, the results indicate that, overall, signal dimension estimations based on SBSS solutions with less kernel functions are more stable. This was already hinted in the former simulations and contrasts with the guidelines presented by \cite{bachoc2018spatial} where it was shown that SBSS solutions based on jointly diagonalizing more than one kernel function yielded more stable solutions.

\begin{table}
\caption{Estimated signal dimension for the moss dataset.}
\centering
\begin{tabular}{cccccc}
  \toprule
Kernel setting & Asym & Perm & Param & Sp Perm & Sp Param  \\ \midrule
Setting 1 & 14 & 14 & 14 & 14 & 14 \\
Setting 2 & 17 & 17 & 18 & 15 & 16 \\
  \bottomrule
\end{tabular}
\label{tab:kola}
\end{table}

\begin{figure}[h]
	\centering
  \includegraphics[scale=0.6]{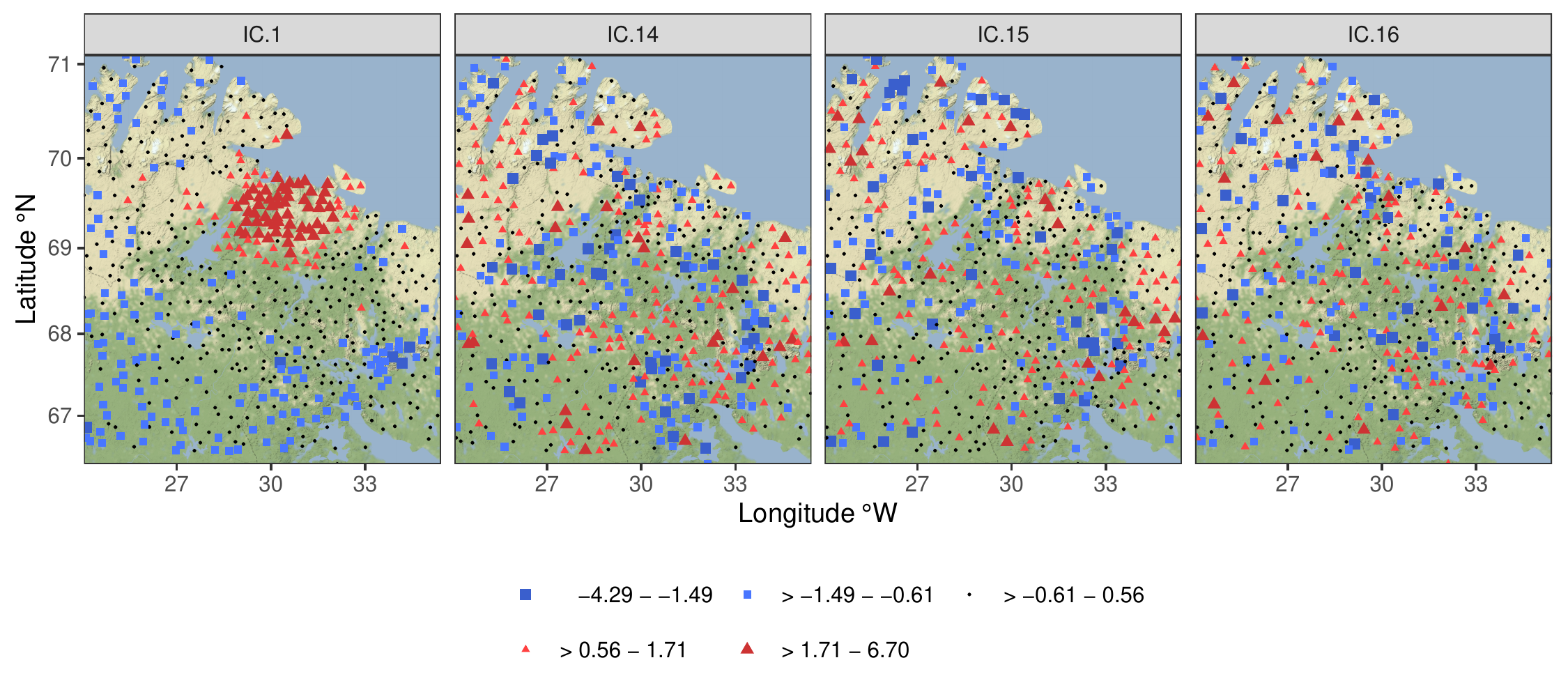}
  \caption{Latent field components of the change-point between estimated signal and noise components for kernel function setting 1. Map tiles by Stamen Design, under CC BY 3.0. Data by OpenStreetMap, under ODbL.}
  \label{fig::kola_1}
\end{figure}

\begin{figure}[h]
	\centering
  \includegraphics[scale=0.6]{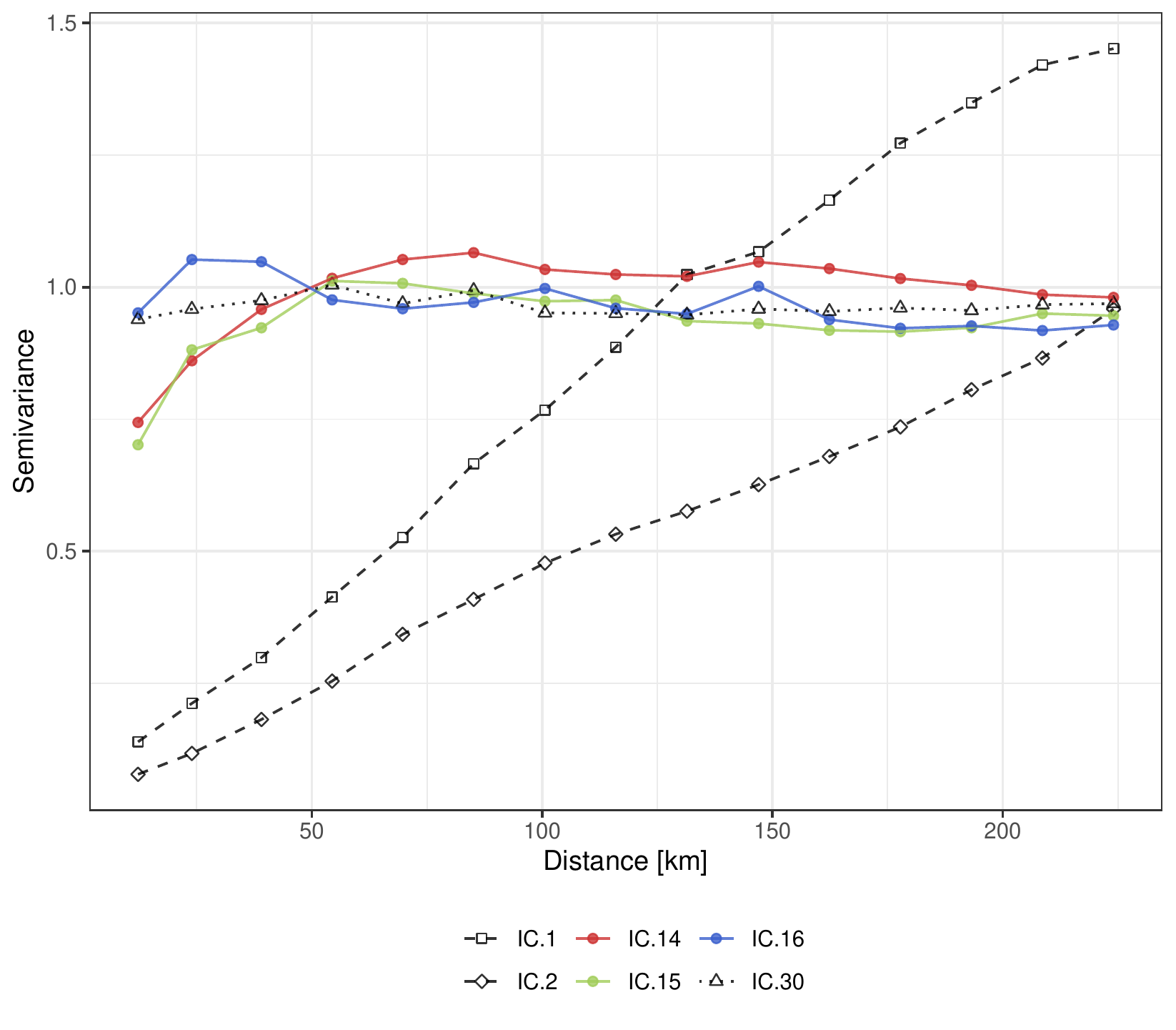}
  \caption{Sample variograms for the latent field components that are first, around the change-point between the estimated signal and noise components and last for kernel function setting 1.}
  \label{fig::kola_vars_1}
\end{figure}

\begin{figure}[h]
	\centering
  \includegraphics[scale=0.6]{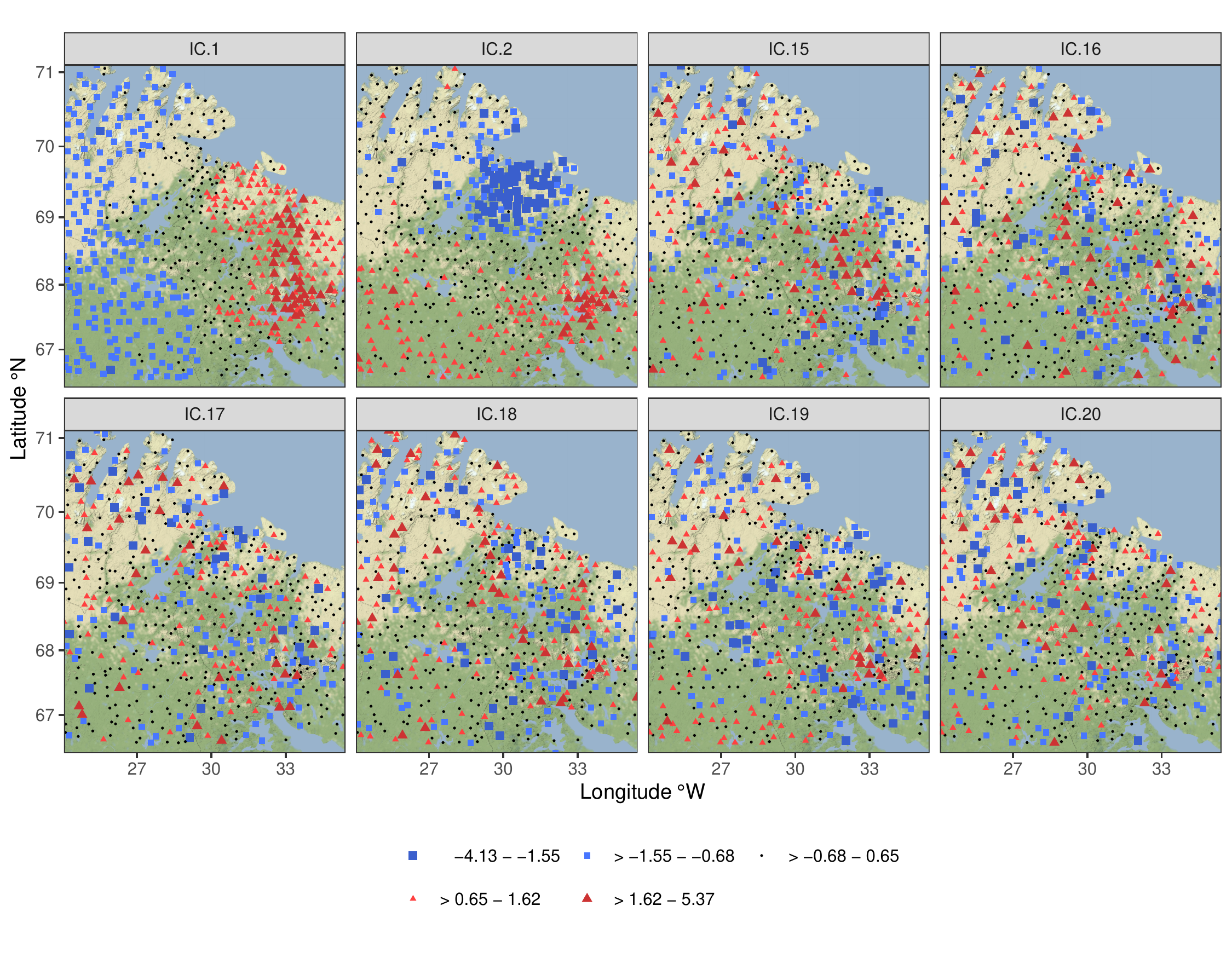}
  \caption{Latent field components of the change-point between the estimated signal and noise components for kernel function setting 2. Map tiles by Stamen Design, under CC BY 3.0. Data by OpenStreetMap, under ODbL.}
  \label{fig::kola_2}
\end{figure}

\begin{figure}[h]
	\centering
  \includegraphics[scale=0.6]{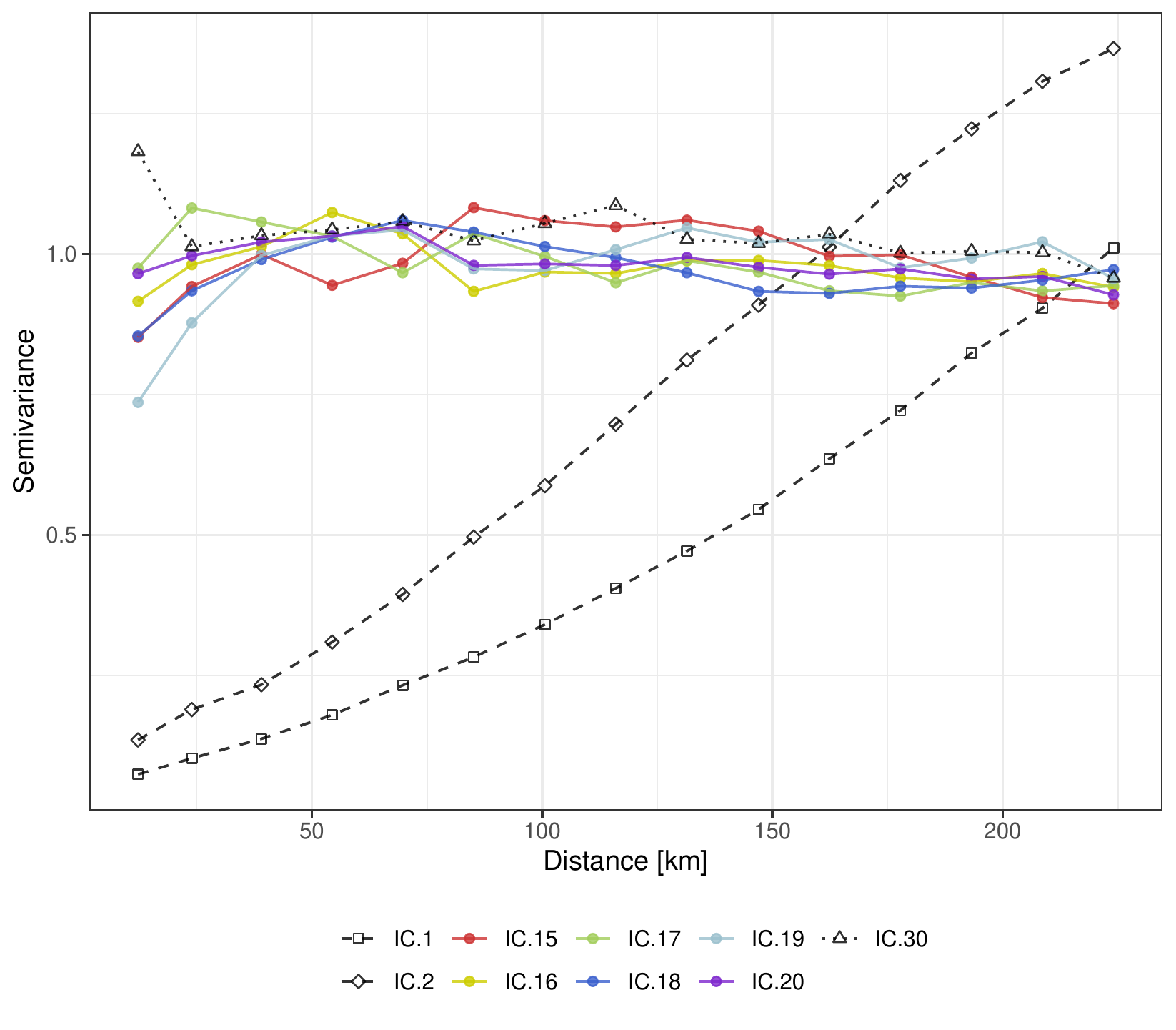}
  \caption{Sample variograms for the latent field components that are first, around the change-point between the estimated signal and noise components and last for kernel function setting 2.}
  \label{fig::kola_vars_2}
\end{figure}

Figures~\ref{fig::kola_1} and \ref{fig::kola_2} show the entries of the latent fields, and Figures~\ref{fig::kola_vars_1} and \ref{fig::kola_vars_2} show the corresponding sample variograms at the change point between the estimated signal and noise for both kernel function settings. Based on visual inspection of Figure~\ref{fig::kola_1}, IC.14-IC.16 might show some weak spatial dependence. However, the sample variograms in Figure~\ref{fig::kola_vars_1} suggest that IC.14 and IC.15 carry a very weak signal, and IC.16 already shows very similar behavior to the last part of the noise, IC.30. Each panel in Figure~\ref{fig::kola_2} could indicate very similar behavior and also illustrates the fact that different signal dimensions are estimated by the different methods. Interestingly, the sample variograms in Figure~\ref{fig::kola_vars_2} show that among IC.15-IC.20, IC.19 had the highest spatial dependence, though all seem to be highly similar to the last components.

\section{Concluding Remarks}\label{Sec:Conclusion}
In this paper, we propose and study testing and estimation methods for the number of latent signal components in the SBSS model. The asymptotic null distributions of the test statistic are given under various conditions without assuming the domain is necessarily regular. A consistent estimator of the dimension based on the sequential tests is also introduced. For small sample cases, different bootstrap strategies are suggested. Besides the theoretical results, the three simulation studies presented in Section \ref{Sec:Sim} demonstrate that our asymptotic tests are comparable to the bootstrap ones in terms of hypothesis testing and estimation. In terms of computation time, our asymptotic method is much faster than the bootstrap ones. When a regular domain structure is used, the computation time can be significantly decreased. 

Our proposed dimension tests in the SBSS context might be very useful for further analysis of the latent fields,  including for various forms of spatial prediction. Indeed, the components of the latent field are uncorrelated, and thus predictions can be carried out on each latent field independently, leading to a reduction from building a single multivariate model to building several univariate ones. This procedure was already investigated and found to be useful by \cite{kriging_paper}. As an additional step, one of our proposed dimension tests can be carried out before the spatial prediction, leading to a reduction of the latent field dimension, which results in the need for even fewer univariate models to be built.

Finally, in further research we plan to develop a ladle estimator (\cite{LuoLi:2016}) for this setting, as well as to develop similar approaches for spatiotemporal data. We will also study the fixed-domain asymptotic properties (\cite{cressie2015statistics}) of SBSS. 

\section{Acknowledgments}
The work of CM, KN, and MY was supported by the Austrian Science Fund P31881-N32.

\appendix

\section{Proofs when $\Om = \mathbf{I}_p$} \label{supplement:section:Om:equal:Ip}

We consider, in this section, the case $\Om = \mathbf{I}_p$.
First, we need to present a useful lemma, which establishes the equivalence between the test statistic and a simpler one that does not involve optimization.

\begin{lemma}\label{lem:eqv}
Under Conditions \ref{c:1}-\ref{c:7}, we have
\[
t_q=t^*_q+o_p(1)\]
as $n\to\infty$, where
\[
t^*_q=\frac{n}{2}\sum_{\ell=1}^k||\wM(f_\ell)_{00}||^2.\]
\end{lemma}

\subsection*{\small \emph{Proof of Lemma \ref{lem:eqv}}}
From Condition \ref{c:7}, there exist $N_0\in\mathbb{N}$ and $F_{inf}>0$ such that for all $n\ge N_0$, for all $\ell=1,\cdots,k, F_{n, f_{\ell}}\ge F_{inf}$. We let $n\ge N_0$ throughout the proof. Assume that 
\[
t_q=t^*_q+o_p(1)\text{ does not hold.}\]
Then there exists $\epsilon>0$ such that
\[
\limsup_{n\to\infty}P(|t_q-t^*_q|\ge\epsilon)\ge\epsilon.\]
We can then extract a subsequence $\phi_1(n)$ such that, along $\phi_1(n)$,
\begin{equation}\label{pf:1}
\liminf_{n\to\infty}P(|t_q-t^*_q|\ge\epsilon)\ge\epsilon.
\end{equation}
The sequence of matrices $\M(f_1), \cdots, \M(f_k)$ are bounded as $n\to\infty$, since for $f\in\{f_1, \cdots, f_k\}$ and $a,b=1,\cdots,p$, from Conditions \ref{c:5} and \ref{c:6},
\begin{align*}
|\M(f)_{a,b}| & \le \frac{1}{n\sqrt{F_{n,f}}}\sum_{i=1}^n\sum_{j=1}^n|f(\s_i-\s_j)|\beta I(a=b)\\
&\le \frac{\beta^2}{\sqrt{F_{n,f}}}\max_{i=1,\cdots,n}\sum_{j=1}^n\frac{1}{1+||\s_i-\s_j||^{d+\alpha}}
\end{align*}
which is bounded from Lemma 4 in \cite{furrer2016asymptotic}. Hence, we can extract a further subsequence $\phi_2(n)$ such that, along $\phi_2(n)$,
\[
\M(f_1)\to \M_{\infty}(f_1), \cdots, \M(f_k)\to \M_{\infty}(f_k).\]
Furthermore, from Condition \ref{c:2}, for $a=1,\cdots,q,$ we have
\[
\sum_{\ell=1}^k \M_{\infty}(f_{\ell})^2_{a,a}>0.\]
We now use some notation from \cite{bachoc2018spatial}. Let $\y$ be the $np\times 1$ vector defined by $\y_{(i-1)p+j}=z_j(\s_i)$, for $i=1,\cdots,n, j=1,\cdots,p.$ For $f\in\{f_0,\cdots,f_k\}$ and for $a,b=1,\cdots,p$, let $\T_{a,b}(f)$ be the $np\times np$ matrix, that we see as a block matrix composed of $n^2$ blocks of sizes $p^2$, and with block $i,j$ equal to
\[
\frac{1}{\sqrt{F_{n,f}}}f(\s_i-\s_j)\frac{1}{2}\{\e_a(p)\e_b(p)^T+\e_b(p)\e_a(p)^T\}.\]
Here $\e_m(p)$ is the $m$-th column basis vector of $\R^p$ for $m=1,\cdots,p$. Then, as observed in the proof of Proposition B.2 in the supplementary material to \cite{bachoc2018spatial}, we have, for $a,b=1,\cdots,p$,
\[
\wM(f)_{a,b}=\frac{1}{n}\y^T\T_{a,b}(f)\y.\]
Furthermore, as observed in the proof of Proposition B.2 in the supplementary material to \cite{bachoc2018spatial}, the largest singular value of $\T_{a,b}(f)$ is bounded as $n\to\infty$. Also, the largest singular value of $\Cov(\y)$ is bounded from Lemma B.1 in the supplementary material to \cite{bachoc2018spatial}. Hence, from Theorem B.1 in the supplementary material to \cite{bachoc2018spatial}, we have
\[
\frac{1}{n}\y^T\T_{a,b}(f)\y-\mathbb{E}\left(\frac{1}{n}\y^T\T_{a,b}(f)\y\right)=\bigO_p\left(\frac{1}{\sqrt{n}}\right).\]
Hence, still along $\phi_2(n)$ we have, for $\ell=1,\cdots,k$,
\begin{equation} \label{eq:hatM:to:M}
\sqrt{n}\left(\wM(f_{\ell})-\M(f_{\ell})\right)=\bigO_p(1).
\end{equation}
Now still along the subsequence $\phi_2(n)$, we apply the proof of Lemma 3 in \cite{virta2018determining}. Although \cite{virta2018determining} address a time series setting, the proof of Lemma 3 there can be directly applied to our setting. Indeed, let us follow the notation of \cite{virta2018determining} and write, for $\ell=1,\cdots,k$,
\begin{equation}\label{eq:Minf}
\begin{pmatrix}
\lambda_{\tau_{\ell}1}^*\\
& \ddots\\
&&\lambda_{\tau_{\ell}q}^*\\
&&& 0\\
&&&&\ddots\\
&&&&&0\end{pmatrix}=\M_{\infty}(f_{\ell}),
\end{equation}
for arbitrary two-by-two distinct $\tau_1,\cdots,\tau_k\in\mathbb{N}$. Let also $\widehat{\U}^T=\wG\wM(f_0)^{1/2}$. With these notations, Assumptions 1 and 3 in \cite{virta2018determining} are satisfied, and the proof of Lemma 3 there can be directly applied. 
Remark that Assumption 3 in \cite{virta2018determining} corresponds to \eqref{eq:hatM:to:M}. 
Remark also that for $\ell = 1,\ldots,k$, $\M(f_{\ell})$ would be fixed and equal to $\M_{\infty}(f_{\ell})$ in the context of \cite{virta2018determining}, while its first $q$ diagonal elements converge to the corresponding elements of $\M_{\infty}(f_{\ell})$ here. Because the last $p-q$ diagonal elements of $\M(f_{\ell})$ are zero, this does not change the arguments in \cite{virta2018determining}. 
The conclusion of the proof of Lemma 3 in \cite{virta2018determining} is then exactly that, along the subsequence $\phi_2(n)$,
\[
t_q=t^*_q+o_p(1).\]
This is in contradiction with \eqref{pf:1}, which concludes the proof.

\subsection*{\small \emph{Proof of Proposition \ref{prop:asymt}}}
Because of Lemma \ref{lem:eqv}, it is sufficient to prove the proposition with $t^*_q$ instead of $t_q$. The statistic $2t^*_q$ is the squared Euclidean norm of the random vector
\[
\sqrt{n}\left(\wM(f_{\ell})_{q+a, q+b}\right)_{\ell=1,\cdots,k,a,b=1,\cdots,p-q}
\]
in dimension $k(p-q)^2$. This random vector is itself equal to
\begin{equation}\label{pf:2}
\sqrt{n}\left(\frac{1}{n\sqrt{F_{n,f_{\ell}}}}\sum_{i,j=1}^nf_{\ell}(\s_i-\s_j)z_{q+a}(\s_i)z_{q+b}(\s_j)\right)_{\ell=1,\cdots,k,a,b=1,\cdots,p-q}.
\end{equation}
With the notation of the proof of Lemma \ref{lem:eqv}, this random vector is equal to
\[
\sqrt{n}\left(\frac{1}{n}\y^T\T_{q+a, q+b}(f_{\ell})\y\right)_{\ell=1,\cdots,k,a,b=1,\cdots,p-q}, \]
where the largest singular value of $\T_{q+a, q+b}(f_{\ell})$ is bounded as $n\to\infty$ for $\ell=1,\cdots,k, a,b=1,\cdots,p-q$ and where the largest singular value of $\Cov(\y)$ is also bounded. Hence, from Theorem B.1 in \cite{bachoc2018spatial}, with $\mathcal{L}_{\M,n}$ the distribution of 
\[
\sqrt{n}\left(\wM(f_{\ell})_{q+a, q+b}-\mathbb{E}\left(\wM(f_{\ell})_{q+a, q+b}\right)\right)_{\ell=1,\cdots,k,a,b=1,\cdots,p-q}\]
we have
\[
d_w(\mathcal{L}_{\M,n}, \mathcal{N}(0, \boldsymbol{\Sigma}_n))\xrightarrow[n\to\infty]{} 0,\]
where $\boldsymbol{\Sigma}_n$ is the covariance matrix of 
$$\sqrt{n}\left(\wM(f_{\ell})_{q+a, q+b}\right)_{\ell=1,\cdots,k,a,b=1,\cdots,p-q}.$$
By the independence of the components of $(z_{q+a}(\s_i))_{a=1,\cdots,p-q,i=1,\cdots,n}$ from Condition \ref{c:1}, and from \eqref{pf:2}, we obtain
\[
\mathbb{E}\left(\wM(f_{\ell})_{q+a, q+b}\right)=0\]
for $\ell=1, \cdots, k, a, b=1,\cdots,p-q$. We also have, for $\ell, \ell'=1, \cdots, k$ and $a, b, a', b' = 1, \cdots, p-q$, from \eqref{pf:2},
\begin{align*}
&n\Cov\left(\wM(f_{\ell})_{q+a, q+b}, \wM(f_{\ell'})_{q+a', q+b'}\right)\\
&=\frac{1}{n\sqrt{F_{n, f_{\ell}}}\sqrt{F_{n, f_{\ell'}}}}\sum_{i,j,i',j'=1}^nf_{\ell}(\s_i-\s_j)f_{\ell'}(\s_{i'}-\s_{j'})\Cov(z_{q+a}(\s_i)z_{q+b}(\s_j), z_{q+a'}(\s_{i'})z_{q+b'}(\s_{j'})).
\end{align*}
From Isserliss' formula applied to the Gaussian vector $(z_{q+a}(\s_i), z_{q+b}(\s_j), z_{q+a'}(\s_{i'}), z_{q+b'}(\s_{j'}))$ with mean vector $\0$, we obtain
\begin{align*}
&\Cov(z_{q+a}(\s_i)z_{q+b}(\s_j), z_{q+a'}(\s_{i'})z_{q+b'}(\s_{j'}))\\
&=\Cov(z_{q+a}(\s_i), z_{q+a'}(\s_{i'}))\Cov(z_{q+b}(\s_j), z_{q+b'}(\s_{j'}))+\Cov(z_{q+a}(\s_i), z_{q+b'}(\s_{j'}))\Cov(z_{q+b}(\s_j), z_{q+a'}(\s_{i'})).
\end{align*}
This yields
\begin{align}
&n\Cov\left(\wM(f_{\ell})_{q+a, q+b}, \wM(f_{\ell'})_{q+a', q+b'}\right) \nonumber \\ 
\label{pf:3}
&=\frac{1}{n\sqrt{F_{n, f_{\ell}}}\sqrt{F_{n, f_{\ell'}}}}\sum_{i,j,i',j'=1}^nf_{\ell}(\s_i-\s_j)f_{\ell'}(\s_{i'}-\s_{j'})\Cov(z_{q+a}(\s_i), z_{q+a'}(\s_{i'}))\Cov(z_{q+b}(\s_j), z_{q+b'}(\s_{j'})) \\ \label{pf:4}
&+\frac{1}{n\sqrt{F_{n, f_{\ell}}}\sqrt{F_{n, f_{\ell'}}}}\sum_{i,j,i',j'=1}^nf_{\ell}(\s_i-\s_j)f_{\ell'}(\s_{i'}-\s_{j'})\Cov(z_{q+a}(\s_i), z_{q+b'}(\s_{j'}))\Cov(z_{q+b}(\s_j), z_{q+a'}(\s_{i'})).
\end{align}
In the sum in \eqref{pf:3}, necessary conditions for a summand to be non-zero are $i=i'$ and $a=a'$ (otherwise the first covariance is zero) and $j=j'$ and $b=b'$ (otherwise the second covariance is zero). In the sum in \eqref{pf:4}, necessary conditions for a summand to be non-zero are $i=j'$ and $a=b'$ (otherwise the first covariance is zero) and $j=i'$ and $b=a'$ (otherwise the second covariance is zero). Hence, we obtain
\begin{align}
&n\Cov\left(\wM(f_{\ell})_{q+a, q+b}, \wM(f_{\ell'})_{q+a', q+b'}\right)\nonumber\\
&=I(a=a')I(b=b')\frac{1}{n\sqrt{F_{n, f_{\ell}}}\sqrt{F_{n, f_{\ell'}}}}\nonumber\\
&\sum_{i,j=1}^nf_{\ell}(\s_i-\s_j)f_{\ell'}(\s_{i}-\s_{j})\Cov(z_{q+a}(\s_i), z_{q+a}(\s_{i}))\Cov(z_{q+b}(\s_j), z_{q+b}(\s_{j}))\nonumber \\ 
&+I(a=b')I(b=a')\frac{1}{n\sqrt{F_{n, f_{\ell}}}\sqrt{F_{n, f_{\ell'}}}}\nonumber\\
&\sum_{i,j=1}^nf_{\ell}(\s_i-\s_j)f_{\ell'}(\s_{j}-\s_{i})\Cov(z_{q+a}(\s_i), z_{q+a}(\s_{i}))\Cov(z_{q+b}(\s_j), z_{q+b}(\s_{j}))\nonumber \\
&=I(a=a')I(b=b')\frac{1}{n\sqrt{F_{n, f_{\ell}}}\sqrt{F_{n, f_{\ell'}}}}\sum_{i,j=1}^nf_{\ell}(\s_i-\s_j)f_{\ell'}(\s_{i}-\s_{j})\nonumber\\
\label{pf:5}
&+I(a=b')I(b=a')\frac{1}{n\sqrt{F_{n, f_{\ell}}}\sqrt{F_{n, f_{\ell'}}}}\sum_{i,j=1}^nf_{\ell}(\s_i-\s_j)f_{\ell'}(\s_{i}-\s_{j})\\
&=I(\ell=\ell')(I(a=a')I(b=b')+I(a=b')I(b=a')),\nonumber
\end{align}
from Condition \ref{c:8}. This means that, for $(\xi_{\ell,a,b})_{\ell=1,\ldots,k,a,b=1,\cdots,p-q}\sim\mathcal{N}(\0, \boldsymbol{\Sigma}_n)$, we have, in distribution,
\[
\sum_{\ell=1}^k\sum_{a,b=1}^{p-q}\xi^2_{\ell,a,b}=\sum_{\ell=1}^k\sum_{a=1}^{p-q}\xi^2_{\ell,a,a}+\sum_{\ell=1}^k\sum_{1\le a<b\le p-q}2\xi^2_{\ell,a,b},\]
where in the two sums on the right-hand-side of the above display, there are $k(p-q)(p-q+1)/2$ independent summands which are squares of Gaussian variables with mean zero and variance 2. Hence, from the continuous mapping theorem, $2t^*_q$ converges in distribution as $n\to\infty$ to 2 times a chi-square-distributed random variable with $k(p-q)(p-q+1)/2$ degrees of freedom. This concludes the proof.

\subsection*{\small \emph{Proof of Proposition \ref{prop:dw}}}
The beginning of the proof is the same as for Proposition \ref{prop:asymt}. We let $\tilde{t}^*_q$ be defined as in Lemma \ref{lem:eqv}, but with $\wM(f)$ replaced by $\widetilde{\M}(f)$ for $f\in\{f_1,\cdots,f_k\}$. Then, the conclusion of this lemma still holds with $(t_q, t^*_q)$ replaced by $(\tilde{t}_q, \tilde{t}^*_q)$. Hence it is sufficient to show the proposition for the distribution $\mathcal{L}_{\tilde{t}^*_q, n}$ of $\tilde{t}^*_q$ instead of the distribution $\mathcal{L}_{\tilde{t}_q, n}$.

The statistic $2\tilde{t}^*_q$ is the squared Euclidean norm of the random vector \[\sqrt{n}\left(\widetilde{\M}(f_{\ell})_{q+a, q+b}\right)_{\ell=1,\cdots,k, a,b=1,\cdots,p-q}\] in dimension $k(p-q)^2$, with mean vector $\0$ and covariance matrix $\widetilde{\Sig}_n$. Let $\mathcal{L}_{\widetilde{\M}, n}$ be the distribution of this vector. As for the proof of Proposition \ref{prop:asymt}, we show that we have
\begin{equation}\label{pf:6}
d_w\left(\mathcal{L}_{\widetilde{\M}, n}, \mathcal{N}(0, \widetilde{\Sig}_n)\right)\xrightarrow[n\to\infty]{} 0.
\end{equation}
Furthermore, similarly as in \eqref{pf:5}, we obtain, for $\ell, \ell'=1,\cdots,k$ and $a,b,a',b'=1,\cdots,p-q$,
\begin{align}
&n\Cov\left(\widetilde{\M}(f_{\ell})_{q+a, q+b}, \widetilde{\M}(f_{\ell'})_{q+a', q+b'}\right)\nonumber\\
&=I(a=a')I(b=b')F_{n,f_{\ell},f_{\ell'}}+I(a=b')I(b=a')F_{n,f_{\ell},f_{\ell'}}\nonumber\\
\label{pf:7}
&=2\Cov(\V_{\ell,a,b}, \V_{\ell',a',b'}).
\end{align}
Hence, if we assume that
\[
d_w\left(\mathcal{L}_{\tilde{t}^*_q, n}, \mathcal{L}_{\V,n}\right)\xrightarrow[n\to\infty]{} 0\text{ does not hold}\]
then there exist $\epsilon>0$ and a subsequence $\phi_1(n)$ such that along $\phi_1(n)$,
\begin{equation}\label{pf:8}
d_w\left(\mathcal{L}_{\tilde{t}^*_q, n}, \mathcal{L}_{\V,n}\right)\ge\epsilon.
\end{equation}
But then we can consider a further subsequence $\phi_2(n)$ along which $F_{n,f_{\ell},f_{\ell'}}$ converges for $\ell,\ell'=1,\cdots,k$ and we can apply the continuous mapping theorem, together with \eqref{pf:6} and \eqref{pf:7}, to obtain a contradiction to \eqref{pf:8}. This concludes the proof.

\subsection*{\small \emph{Proof of Corollary \ref{coro:dw}}}
Consider i.i.d. standard Gaussian variables $(X_{\ell,a,b})_{\ell=1,\cdots,k,1\le a\le b\le p-q}$, independent of $\V$. We let for $\ell=1,\cdots,k$ and $1\le a \le b\le p-q$,
\[
W^{(\ell)}_{a,b}=\begin{cases}
X_{\ell,a,b} & \text{if } F_{n,f_{\ell}}=0\\
\frac{1}{\sqrt{F_{n,f_{\ell}}}}\V_{\ell,a,a} & \text{if } F_{n,f_{\ell}}\neq 0\text{ and } a=b\\
\frac{1}{\sqrt{2 F_{n,f_{\ell}}}}(\V_{\ell,a,b}+\V_{\ell,b,a}) & \text{if } F_{n,f_{\ell}}\neq 0\text{ and } a\neq b
\end{cases}.\]
Then, it can be checked that the $(W^{(\ell)}_{a,b})_{\ell=1,\cdots,k,1\le a\le b\le p-q}$ are i.i.d. standard Gaussian variables. Furthermore,
\[
\sum_{\ell=1}^k\sum_{a,b=1}^{p-q}\V^2_{\ell,a,b}=\sum_{\ell=1}^k\sum_{a=1}^{p-q}\V^2_{\ell,a,a}+\sum_{\ell=1}^k\sum_{1\le a<b\le p-q}2\V^2_{\ell,a,b},\]
since $\V_{\ell,a,b}$ has correlation one with $\V_{\ell,b,a}$ and same variance. We then have
\begin{align*}
\sum_{\ell=1}^k\sum_{a,b=1}^{p-q}\V^2_{\ell,a,b}
&=
\sum_{\ell=1}^k\sum_{a=1}^{p-q}I(F_{n,f_{\ell}}\neq 0)\V^2_{\ell,a,a}+\sum_{\ell=1}^k\sum_{1\le a<b\le p-q}I(F_{n,f_{\ell}}\neq 0)(2\V^2_{\ell,a,b}) 
\\
 &=
 \sum_{\ell=1}^k\sum_{a=1}^{p-q}I(F_{n,f_{\ell}}\neq 0)\V^2_{\ell,a,a}+\sum_{\ell=1}^k\sum_{1\le a<b\le p-q}I(F_{n,f_{\ell}}\neq 0)\left(\frac{1}{2}(\V_{\ell,a,b}+\V_{\ell,b,a})^2\right)
 \\
&=
\sum_{\ell=1}^k\sum_{a=1}^{p-q}\left(\sqrt{F_{n,f_{\ell}}}W_{a,a}^{(\ell)}\right)^2+\sum_{\ell=1}^k\sum_{1\le a<b\le p-q}\left(\sqrt{F_{n,f_{\ell}}}W_{a,b}^{(\ell)}\right)^2.
\end{align*}
This concludes the proof.

\subsection*{\small \emph{Proof of Proposition \ref{prop:power}}}
Clearly, the test statistics are monotonous, in the sense that $t_r\ge t_{r'}$ for $r\le r'$. Indeed, $t_r$ is obtained from the norm of a matrix and $t_{r'}$ is obtained from the norm of a submatrix of this matrix.

Hence, from Proposition \ref{prop:asymt}, the conclusion of Proposition \ref{prop:power} holds for $r\ge q$. Let now $r<q$. If $q=0$, there is nothing that needs to be proved. If $q\ge 1$, from the previously remarked monotonicity, it is sufficient to consider $r=q-1$. Let, in view of Condition \ref{c:2}, 
\[
b = \frac{1}{2}\liminf_{n \to \infty}\min_{a=1,\ldots,q}\sum_{\ell = 1}^k\M(f_{\ell})_{a,a}^2 >0.
\]
Assume that $ t_r/n \geq b + o_p(1)$ does not hold. Then there exists $c < b$ such that $\mathbb{P}(  t_r/n \geq c )$ does not go to one as $n \to \infty$.
 Then there exists a subsequence $\phi(n)$ and $c < b$ such that, along $\phi(n)$,
\begin{equation} \label{eq:in:proof:power:to:be:contradicted}
\limsup\mathbb{P}(   t_r/n \geq c )< 1.
\end{equation}
As in the proof of Lemma \ref{lem:eqv}, let us extract a further subsequence $\phi_2(n)$ such that, along $\phi_2(n)$, 
\[
\M(f_1) \to \M_{\infty}(f_1), \ldots, \M(f_k) \to \M_{\infty}(f_k).
\]
We have 
\begin{equation} \label{eq:in:proof:power:larger:two:b}
\sum_{\ell = 1}^k\M_{\infty}(f_{\ell})_{q,q}^2 \geq 2 b .
\end{equation}
We also have, along $\phi_2(n)$,
\begin{equation} \label{eq:in:proof:power:tr:larger}
\frac{t_r}{n}\geq \frac{1}{2}\sum_{\ell=1}^k\widehat{\D}_{\ell,q,q}^2.
\end{equation}
Similarly to the proof of Proposition \ref{prop:asymt}, the proof of Lemma 5 in \cite{virta2018determining} can be directly applied to our setting, with \eqref{eq:Minf}. In the proof of this lemma, it is shown that, along $\phi_2(n)$,
\begin{equation} \label{eq:in:proof:power:D:to:M}
\sum_{\ell=1}^k\widehat{\D}_{\ell,q,q}^2\xrightarrow[n \to \infty]{p}\sum_{\ell = 1}^k\M_{\infty}(f_{\ell})_{q,q}^2.
\end{equation}
Hence from \eqref{eq:in:proof:power:larger:two:b}, \eqref{eq:in:proof:power:tr:larger} and \eqref{eq:in:proof:power:D:to:M}, we obtain, along $\phi_2(n)$,
\[
\frac{t_r}{n}\geq b + o_p(1),
\]
which is in contradiction with\eqref{eq:in:proof:power:to:be:contradicted}. Hence $ t_r/n \geq b + o_p(1)$ which concludes the proof.

\subsection*{\small \emph{Proof of Proposition \ref{prop:estimating:q}}}
Because Proposition \ref{prop:power} has been established, the proof is the same as in \cite{virta2018determining}. We have, for $r = 1 , \ldots , q-1$,
\[
\mathbb{P} \left(  t_r \leq c_n \right) \xrightarrow[]{n \to \infty} 0,
\]
from Proposition \ref{prop:power}. Furthermore, 
\[
\mathbb{P} \left(  t_q \leq c_n \right) \xrightarrow[]{n \to \infty} 1,
\]
from Proposition \ref{prop:power}. This shows that $\hat{q} = q$ with probability going to one which concludes the proof.

\subsection*{\small \emph{Proof of Proposition \ref{prop:nonzerom}}}
From the expression of $\widebar{\M}(f)$, we can assume without loss of generality that $\mu_1=\cdots=\mu_p=0$. Then, from Lemma B.8 in the supplementary material to \cite{bachoc2018spatial}, we have for $f\in\{f_0,f_1,\cdots,f_k\}$,
\[
\widebar{\M}(f)-\wM(f)=\bigO_p\left(\frac{1}{n}\right).\]
Hence the proofs of Propositions \ref{prop:asymt}, \ref{prop:dw}, \ref{prop:power} and \ref{prop:estimating:q}
and of Corollary \ref{coro:dw} are carried out with $\widebar{\M}(f)$ for $f\in\{f_0,f_1,\cdots,f_k\}$, exactly as with $\wM(f)$ for $f\in\{f_0,f_1,\cdots,f_k\}$.

\section{Proofs for a General $\Om$} \label{supplement:section:general:Om}

Let, for $f \in \{f_0,f_1,\ldots,f_k\}$, $\widehat{\M}_{\boldsymbol{I}_p}(f)$ be defined as in \eqref{eq:slc}, but in the case where $\Om = \boldsymbol{I}_p$ (equivalently with $\x$ replaced by $\z$). 
Let
$\widehat{\M}_{\Om}(f)$ be defined as in \eqref{eq:slc} for $f \in \{f_0,f_1,\ldots,f_k\}$ (insisting on the dependency in $\Om$ in the notation).
We then have $\widehat{\M}_{\Om}(f)  = \Om \widehat{\M}_{\boldsymbol{I}_p}(f) \Om^T $.

Let
\[
t_{\Om}=\frac{n}{2}\sum_{ \ell=1}^k||\wD_{\Om,\ell,00}||^2\]
be the test statistic in \eqref{eq:t}, insisting on the dependency in $\Om$ in the notation. Recall that for $\ell=1,\cdots,k, \wD_{\Om,\ell,00}$ is the lower $(p-q)\times(p-q)$ diagonal block of 
\[
\wD_{\Om,\ell}=\wG_{\Om}\wM_{\Om}(f_{\ell})\wG_{\Om}^T,\]
with
\[
\widehat{\boldsymbol \Gamma}_{\Om}
\in 
\argmax_{
\substack{
 \boldsymbol{\Gamma}_{\Om}: \boldsymbol{\Gamma}_{\Om} \widehat{\M}_{\Om}(f_0) \boldsymbol {\Gamma}_{\Om}^\top = \I_p  \\
\boldsymbol{\Gamma}_{\Om} \mbox{\small \; has rows } \boldsymbol{\gamma}_{\Om,1}^\top,\ldots,\boldsymbol\gamma_{\Om,p}^\top \\
\left( \sum_{\ell=1}^k \{ \boldsymbol{\gamma}_{\Om,j}^\top \widehat{\M}_{\Om}(f_\ell) \boldsymbol{\gamma}_{\Om,j} \}^2
\right)_{j = 1,\ldots , p}
~ \mbox{are in descending order}
} 
 }
 \sum_{\ell=1}^k \sum_{j=1}^p \{ \boldsymbol{\gamma}_{\Om,j}^\top \widehat{\M}_{\Om}(f_\ell) \boldsymbol{\gamma}_{\Om,j} \}^2.
\]

We then also have
\begin{flalign*}
&
\widehat{\boldsymbol \Gamma}_{\Om} \Om
\in 
\argmax_{
\substack{
 \boldsymbol{\Gamma}_{\Om} \Om: \boldsymbol{\Gamma}_{\Om}
 \Om  \widehat{\M}_{\boldsymbol{I}_p}(f_0) \left( \boldsymbol {\Gamma}_{\Om} \Om \right)^T = \I_p  \\
\boldsymbol{\Gamma}_{\Om} \Om \mbox{\small \; has rows } \boldsymbol{\gamma}_{\Om,1}^T \Om,\ldots,\boldsymbol\gamma_{\Om,p}^T \Om \\
\left( \sum_{\ell=1}^k \{ 
\left( \boldsymbol{\gamma}_{\Om,j}^\top \Om \right) \widehat{\M}_{\boldsymbol{I}_p}(f_\ell) \left( \boldsymbol{\gamma}_{\Om,j}^T \Om \right)^T \}^2
\right)_{j = 1,\ldots , p}
~ \mbox{are} \\
\mbox{in descending order}
} 
 }
 \sum_{\ell=1}^k \sum_{j=1}^p \{ 
\left( \boldsymbol{\gamma}_{\Om,j}^T \Om \right) \widehat{\M}_{\boldsymbol{I}_p}(f_{\ell}) \left( \boldsymbol{\gamma}_{\Om,j}^T \Om \right)^T \}^2.
\end{flalign*}

Furthermore, we have $t_{\Om}=t_{\I_p}$ where
\[
t_{\I_p}=\frac{n}{2}\sum_{\ell=1}^k||\wD_{\I_p,\ell,00}||^2\]
where, for $\ell=1,\cdots,k,\wD_{\I_p,\ell,00}$ is composed by the lower $(p-q)\times(p-q)$ diagonal block of the matrix
\[
\wD_{\I_p,\ell}=\left(\wG_{\Om}\Om\right)\wM_{\I_p}(f_{\ell})\left(\wG_{\Om}\Om\right)^T.\]
The test statistic $t_{\I_p}$ has the asymptotic distribution given in Propositions \ref{prop:asymt} and \ref{prop:dw} and in Corollary \ref{coro:dw}, from the proofs of these propositions and corollary in the case where $\Om=\I_p$. This concludes the proofs of Propositions \ref{prop:asymt} and \ref{prop:dw} and of Corollary \ref{coro:dw} for a general $\Om$. 
This also directly extends the proofs of Propositions \ref{prop:power} and \ref{prop:estimating:q} to the case of a general $\Om$. Finally, the proof of Proposition \ref{prop:nonzerom} for a general $\Om$ is
obtained similarly as above, by observing that, for $f \in \{ f_0,f_1,\ldots,f_k\}$, $\widebar{\M}_{\Om}(f) = \Om \widebar{\M}_{\I_p}(f) \Om^\top$, by extending the notation $\widehat{\M}_{\Om}(f)$ and $\widehat{\M}_{\I_p}(f)$ above to $\widebar{\M}_{\Om}(f)$ and $\widebar{\M}_{\I_p}(f)$.

\bibliographystyle{agsm}



\end{document}